\documentclass[11pt]{article}

\input{amssym.def}
\input{amssym}

\def\theequation{\thesection.\arabic{equation}}
\makeatletter
\@addtoreset{equation}{section}
\makeatother

\def\be{\begin{equation}}
\def\ee{\end{equation}}
\def\ba{\begin{eqnarray}}
\def\ea{\end{eqnarray}}
\def\lb{\label}

\def\tr{{\rm Tr}}
\def\Tr#1{{\rm Tr\str{-1.3}}_{R^{\mbox{\scriptsize$(#1)$}}}}
\def\str#1{\rule[#1mm]{0pt}{1mm}}

\def\tbl#1#2{{\ifmmode \left\{\!\!
\begin{array}{c}
\scriptstyle #1\\[-2pt] \raisebox{2pt}{$\scriptstyle #2$}
\end{array}
\!\!\right\}
\else
\raisebox{2pt} {\scriptsize$\left\{\!\!\!
\begin{array}{c}
#1\\[-1pt] #2
\end{array}
\!\!\!\right\}$}\fi}}

\def\schf#1{s_{\raisebox{-0.08mm}{$_{#1}$}}}

\newcounter{theorem}
\newtheorem{pred}[theorem]{Proposition}
\newtheorem{lem}[theorem]{Lemma}
\newtheorem{opred}[theorem]{Definition}
\newtheorem{teor}[theorem]{Theorem}
\newtheorem{remark}[theorem]{Remark}
\newtheorem{cor}[theorem]{Corollary}

\begin{document}

\title{The $GL(m|n)$ type quantum matrix algebras II:\\
the structure of the characteristic subalgebra\\
and its spectral parameterization}

\author{
\rule{0pt}{7mm} Dimitri Gurevich\thanks{gurevich@univ-valenciennes.fr}\\
{\small\it USTV, Universit\'e de Valenciennes,
59304 Valenciennes, France}\\
\rule{0pt}{7mm} Pavel Pyatov\thanks{pyatov@thsun1.jinr.ru}\\
{\small\it Max Planck Institute for Mathematics,
Vivatsgasse 7, D-53111 Bonn, Germany}\\[-1pt]
{\small\it and}\\[-1pt]
{\small\it Bogoliubov Laboratory of Theoretical Physics,
JINR, 141980 Dubna, Moscow region, Russia}\\
\rule{0pt}{7mm} Pavel Saponov\thanks{Pavel.Saponov@ihep.ru}\\
{\small\it Division of Theoretical Physics, IHEP, 142284
Protvino, Russia}
}

\date{}

\maketitle

\begin{abstract}
In our previous paper \cite{GPS2} the Cayley-Hamilton
identity for the $GL(m|n)$ type quantum matrix algebra was obtained.
Here we continue investigation of that identity.
We derive it in three alternative forms and, most importantly,
we obtain it in a factorized form. The factorization leads to a
separation of the spectra of the quantum supermatrix into the
"even" and "odd" parts. The latter, in turn, allows us to
parameterize the characteristic subalgebra
(which can also be called the subalgebra of spectral invariants)
in terms  of the supersymmetric polynomials in the eigenvalues
of the quantum supermatrix.
For our derivation we use two auxiliary results which
may be of independent interest. First, we calculate
the multiplication rule for the linear basis
of the Schur functions $s_\lambda(M)$ for the
characteristic subalgebra of the Hecke type quantum matrix algebra.
The structure constants in this basis are the Littlewood-Richardson
coefficients.
Second, we derive a series of bilinear relations in the graded ring
$\Lambda$ of Schur symmetric functions in countably many variables
(see \cite{Mac}).
\end{abstract}

\newpage

\tableofcontents

\vskip 15mm
\section{Introduction}
\lb{sec1}

In the present paper, we continue the investigation of
the supersymmetric $GL(m|n)$ type quantum matrix (QM) algebras
initiated in \cite{GPS2}. Let us recall briefly the history of the subject.

The first examples of the QM~algebras
were considered in the seminal papers of V. Drinfel'd \cite{D} and
L. Faddeev, N. Reshetikhin
and L. Takhtajan \cite{FRT}. There, a particular family of QM~algebras ---
the algebras
of quantized functions on the groups, shortly called the RTT~algebras,
were defined.
Soon after, another important subclass of QM~algebras ---
the  reflection equation (RE~) algebras, were introduced into consideration
(see, e.g., \cite{KS,KSas}).
The general definition of the QM~algebras
was found by L.Hlavaty who aimed at giving a unified
description for RTT~ and RE~algebras \cite{Hl}.
This idea might seem quite strange at  first glance (the
representation theories of the RTT~ and the RE~algebras are very different).
At the same time, the structure investigations carried out separately
for the RTT \cite{EOW,Zh,IOPS} and  the RE~algebras
\cite{NT,PS,GPS1} reveal a remarkable similarity of  both the algebras
to the classical matrix algebra.
Namely, it turns out that the RE~ and the RTT~
families admit a noncommutative generalization
of the Cayley-Hamilton theorem and for the matrices of generators
in both the cases
a noncommutative analogue of their spectra can be constructed.
Having this in mind, the general definition of the QM~algebras was
independently reproduced in \cite{IOP2}
and the noncommutative version of the Cayley-Hamilton theorem
was derived for the QM~algebras of the general linear type
(see \cite{IOP1,IOP2,IOP3}).

The family of $GL(m)$ type QM~algebras was a good
case to start with. An investigation of the other classical
series of the QM~algebras falls into two cases --- the case of the Hecke type
QM~algebras and the case of the Birman-Murakami-Wenzl (BMW)
type QM~algebras. The  difference is in the
choice of a quotient of the group algebra of the braid group
which enters (through its R-matrix representation) into
the QM~algebra definition. The Hecke case contains the
general linear type and its supersymmetric generalization --- the $GL(m|n)$
type QM~algebras. The BMW case includes orthogonal- and symplectic- type
QM~algebras and their supersymmetric analogues. An
investigation of the BMW case
was started in \cite{OP2}, where the Cayley-Hamilton identity
and the spectra
of the orthogonal- and symplectic- type QM~algebras were identified.

The supersymmetric $GL(m|n)$ type QM algebra was studied in our previous
paper \cite{GPS2}. In  that paper, we gave a proper definition of the family of
the $GL(m|n)$ type QM~algebras and proved
the Cayley-Hamilton identity for them.
Our work may be viewed as a generalization of both
the results by I. Kantor and I. Trishin on the Cayley-Hamilton
equation for the supermatrices \cite{KT1,KT2} (the invariant Cayley-Hamilton
equation in their terminology), and by P.D. Jarvis and
H.S. Green on the characteristic identities
for the general linear Lie superalgebras \cite{JG}.

Still lacking in the $GL(m|n)$ case is the identification
of the spectrum of the quantum supermatrices.\footnote{
Here we put the problem for the generic QM~algebra. For the subfamily of
RE~algebras and at the level of finite dimensional representations it was
considered in \cite{GL}.}
Alternatively, one can ask for a proper parameterization of the characteristic
subalgebra of the $GL(m|n)$ type QM~algebra (the abelian subalgebra
of the QM~algebra which
the coefficients of the Cayley-Hamilton identity belong to).
This problem is addressed in the present work. First, we
investigate in detail the structure of the characteristic
subalgebra in the Hecke case. Then, we derive a series of
bilinear relations in the graded ring $\Lambda$
of Schur symmetric functions
in countably many variables (for the definition see \cite{Mac}).
These combinatorial relations may be of independent interest.

The structure of the paper is as follows.
In the next section,  subsection
\ref{sec2.1}, we derive the multiplication rule for the set of
linear basic elements of the Hecke type characteristic subalgebra ---
the so-called Schur functions $s_\lambda(M)$ (the notation is explained below).
The structure constants in this basis are just the Littlewood-Richardson
coefficients. In other words, we define the homomorphic map
from the ring of symmetric functions $\Lambda$ to the characteristic
subalgebra of the Hecke type QM~algebra. To efficiently apply this map
in the $GL(m|n)$ case, we need a series of bilinear relations for
the Schur symmetric functions $s_\lambda\in \Lambda$.\footnote{
There should be no confusion between the elements $s_\lambda\in\Lambda$
and their homomorphic images $s_\lambda(M)$ in the characteristic
subalgebra. The argument in the latter notation
is used for distinguishing purposes. It
refers to the matrix of generators of the QM~algebra.
}
They are proved in  subsection \ref{sec2.2}.
For  derivation we use the Jacobi-Trudi formulas for the Schur functions
$s_\lambda$ and apply the Pl{\"u}cker relations. The same method was used
in \cite{LWZ,Kl} for the derivation of  different
bilinear relations for the Schur functions.
We also remark that our bilinear relations certainly have
common roots with the factorization formula
for the supersymmetric functions \cite{BR,PrT}.

In section 3, we derive three alternative expressions for the Cayley-Hamilton
identity for the $GL(m|n)$ type QM~algebra.
In  subsection \ref{sec3.1}, the bilinear identities of  subsection
\ref{sec2.2} are used to factorize the $GL(m|n)$ type characteristic
identity into a product of two terms.  Let us stress that the
factorization is achieved without extending the algebra by the
eigenvalues of the quantum supermatrix.
To the best of our knowledge, this fact has not been observed before even in the
classical supermatrix case.
The factorization allows us to separate "even" and "odd"
eigenvalues of the quantum supermatrix in a covariant manner.
That is, we do not specify explicitly the ${\Bbb Z}_2$-grading
for the components of the quantum supermatrix. Instead, we observe the
"manifestation of even and odd variables"
in the factorization property of the characteristic polynomial.
Two more versions of the Cayley-Hamilton theorem are
presented in  subsection \ref{sec3.2}.
They are given in terms of the (skew-)symmetric powers of the
quantum matrices\footnote{
The notion of the skew-symmetric power of the matrix
was suggested by A.M. Lopshits (see \cite{GGB}, p.342.)}
and generalize the corresponding results of \cite{KT2,T}
to the case $q\neq 1$.
Yet another series of bilinear relations for the
Schur symmetric functions $s_\lambda$
is used here for derivations (see lemma \ref{lem:4}).
These relations are also applied in the last section
for parameterization of the Schur functions
$s_\lambda(M)$.

Finally, in section \ref{sec4}, we compute expressions
for the coefficients of the $GL(m|n)$ type
Cayley-Hamilton identity
in terms of the quantum matrix eigenvalues.
The resulting  parameterization
is given in terms of the supersymmetric polynomials \cite{Stem}
(see also \cite{Mac}, section 1.3, exercises 23 and 24).
It is worth mentioning that the supersymmetric polynomials were originally
introduced
by F. Berezin  \cite{Ber1,Ber2}
for a description of invariant polynomials
on the Lie superalgebra ${\frak g\frak l}(m|n)$
(see also \cite{Ser1} and references therein).

Some auxiliary q-combinatorial formulae which we need for derivations  in
section \ref{sec2.1} are proved in the appendix.

\smallskip
Throughout this text we keep the notation of the paper \cite{GPS2}.
When referring to formulae from that paper we use the shorthand
quotation, e.g., symbol (I-3.21) refers to formula (21) from  section 3
of \cite{GPS2}.
For reader's convenience  in the rest of the introduction
we collect a list of notation, definitions  and results mainly  from
\cite{GPS2}.

\smallskip
Let $V$ be a finite dimensional $\Bbb C$-linear
space, $\dim V = N$.
Consider a pair of elements $R, F\in {\rm Aut}(V^{\otimes 2})$.
Fixing some basis $\{v_i\}_{i=1}^N$ in the space $V$ we identify
operators $R$ and $F$ with their matrices in that basis.
We use the shorthand matrix notation of \cite{FRT}. I.e., we write
$R_i$ (or, sometimes, more explicitly $R_{i\, i+1}$) for the matrix of the
operator
${\rm Id}^{\otimes (i-1)}\otimes R\otimes {\rm Id}^{\otimes (k-i-1)}$
acting in the space $V^{\otimes k}$.
Here $\rm Id\in {\rm Aut(V)}$ denotes the identity operator.
The integer $k$
is not shown in the matrix notation. In each
particular formula the actual value of $k$ can be easily reconstructed.
Few more conventions: $I$ is the identity matrix;
$P\in {\rm Aut}(V^{\otimes 2})$ is the permutation automorphism
($P(u\otimes v) = v\otimes u$).

The pair of operators $R$ and $F$ can be used as an initial data set for the
QM~algebra, provided they satisfy the following conditions
\begin{enumerate}
\item[i)]
The matrices of both operators $R$ and $F$ are {\em strict skew invertible}.
The skew invertibility means, say for $R$, the existence of an operator
$\Psi^R\in {\rm End}(V^{\otimes 2})$ such that
$
{\rm Tr}_{(2)}R_{12}\Psi_{23}^R  = P_{13}\, ,
$
where the subscript in the notation of the trace shows the number of the
space  $V$,  where the trace is evaluated
(here we adopt labelling  $V^{\otimes k} := V_1\otimes V_2\otimes\dots \otimes
V_k$).
The strictness condition implies invertibility of
an element $D_1^R:={\rm Tr}_{(2)}\Psi^R_{12}$:
$D^R\in {\rm Aut}(V)$.

With the matrix $D^R$ one then defines the
{\em R-trace} operation $\tr_R : \; {\rm Mat}_N(W)\rightarrow W$
$$
\tr_R(X) := \sum_{i,j=1}^N{D^R}_i^jX_j^i,\quad  X\in
{\rm Mat}_N(W),
$$
where $W$ is any linear space (in considerations below $W$ is the space
of the QM~algebra).

\item[ii)]
The operators  $R$ and $F$ are the {\em R-matrices}, that means they satisfy
the {\em Yang-Baxter equations}
$$
R_1R_2R_1 = R_2R_1R_2\, ,\qquad F_1F_2F_1 = F_2F_1F_2\, .
$$
\item[iii)]
The operators  $R$ and $F$ form a {\em compatible pair} $\{R,F\}$
(the order of operators in this notation is essential)
$$
R_1F_2F_1 = F_2F_1R_2\, ,
\qquad
F_1F_2R_1 = R_2F_1F_2\, .
$$
\end{enumerate}

Given the pair $\{R,F\}$ satisfying conditions i)--iii)
the {\em quantum matrix  algebra} ${\cal M}(R,F)$
is defined as a unital  associative  algebra which is generated by $N^2$
components of the
matrix  $\|M^i_j\|_{i=1}^N$ subject to the relations\footnote{
In \cite{GPS2} we also demand skew invertibility of an operator
$R_f := F^{-1} R^{-1} F$ in the definition of the QM~algebra.
As is proved in \cite{OP2} (see lemma 3.6),
the latter condition  is a consequence of
i)--iii).}
\be
R_1M_{\overline 1}M_{\overline 2} = M_{\overline 1}
M_{\overline 2}
R_1\, .
\ee
Here we used the iterative procedure
\be
\lb{RMM}
M_{\overline 1} = M, \quad M_{\overline {k+1}}  =
F_k M_{\overline k} F^{-1}_k
\ee
for the production of  copies $M_{\overline{k}}$
of the matrix $M$.
The defining relations (\ref{RMM}) then imply the same type relations
for any consecutive pair of the copies of $M$ (see  lemma I-4)
\be
R_k\,M_{\overline k}M_{\overline{k+1}} =
M_{\overline k}M_{\overline{k+1}}\, R_k.
\label{mk-r}
\ee

Imposing additional conditions on the R-matrix $R$
we then extract specific series of the QM~algebras.
\begin{enumerate}
\item[iv)]
Demanding $R$ to be the {\em Hecke type}
R-matrix, that means its minimal polynomial to be of the second order
\be
\lb{Hecke}
(R+q^{-1}I)(R-qI) = 0\, , \quad q\in \{{\Bbb C}\setminus 0\}\, ,
\ee
we specify to the {\em Hecke type} QM~algebra.
The $\Bbb C$-number $q$  becomes the parameter of the
algebra.
\item[v)]
Given a Hecke type R-matrix (\ref{Hecke}), one can construct a series of
{\em R-matrix representations} of the {\em Hecke algebras\footnote{
A brief description of the Hecke algebras,
their R-matrix representations,
the primitive idempotents
and the basis of matrix units
is given in \cite{GPS2}, sections 2 and 3.
For a more detailed exposition of the subject the reader is referred
to \cite{R,OP1} and to references therein.
}} ${\cal H}_p(q)$
\be
\lb{rhoR}
\rho_R: {\cal H}_p(q)\rightarrow {\rm End}
(V^{\otimes p}), \quad p=2,3,\dots\,.
\ee
Let us impose an additional restriction on the parameter $q$
\be
\lb{semisimple}
q^{2k}\neq 1, \qquad  k=2,3,\dots\, ,
\ee
which ensures the algebras ${\cal H}_p(q)$, $p=2,3,\dots\,$,
to be semisimple.
Then we can further specify to a series of the {\em $GL(m|n)$ type}
QM~algebras. For their definition we use a set of the primitive idempotents
$E^\lambda_\alpha\in H_p(q)$
labelled by the standard Young tableaux \tbl{\lambda}{\alpha}, where
$\lambda\vdash p$ is a partition of $p$, and
index $\alpha$ enumerates different
standard  tableaux corresponding to the partition $\lambda$
(see section I-2).
The {\em $GL(m|n)$ type} QM~algebra
is characterized by the following conditions
\begin{enumerate}
\item[a)]
the representations $\rho_R$  (\ref{rhoR}) are faithful for all
$p < (m+1)(n+1)$;
\item[b)]
for $p\geq (m+1)(n+1)$ the kernel of $\rho_R$
is generated by (any one of) the primitive
idempotents $E^{((n+1)^{(m+1)})}_\alpha$
corresponding to the rectangular Young diagram
$((n+1)^{(m+1)})$;
\item[c)]
the Schur function $s_{(n^m)}(M)$ (see definition below)
corresponding to the rectangular Young diagram $(n^m)$ is an
invertible element of the QM~algebra.\footnote{
This condition was not imposed in \cite{GPS2}.
We will need it now for the spectral parameterization of
the characteristic subalgebra (see eqs.(\ref{def:mu}),
(\ref{def:nu})).}
\end{enumerate}
\end{enumerate}

Two comments are in order:
\begin{itemize}
\item
It is well known that the classical ($q=1$) $GL(m|n)$ type supermatrices
satisfy properties v)-a) and v)-b)
(for a proof see \cite{Ser2}, theorem 2, and \cite{KT1}, theorem 1.2).
Hence, using the deformation arguments
we can make sure that
these properties remain valid for a variety of the QM-algebras
related to the {\em standard} $GL(m|n)$ type R-matrices described, e.g., in
\cite{DKS,I}.
For our applications (at least) it is convenient
to use v)-a) and v)-b) as defining conditions
for the $GL(m|n)$ type QM-algebras.
It seems verisimilar that (for  $q$ generic)
any Hecke type QM-algebra is of $GL(m|n)$ type (for some values of integers
$m$ and $n$). We are going  to further argue this point in a separate work.
\item
Notice that relation $m+n=N\equiv\dim V$
between the algebra parameters $m$, $n$ and $N$
is not assumed in the definition.
Although it is indeed satisfied in many
examples (say, for the QM~algebras constructed by the standard $GL(m|n)$ type
R-matrices), there are known exceptions from this rule.
A series of counter-examples was constructed in \cite{Gur}.
\end{itemize}

{}From now on we restrict ourselves to considering
the Hecke type QM~algebras with the parameter $q$
satisfying condition (\ref{semisimple}).

The {\em characteristic subalgebra}
${\rm Char}(R,F)$
of the QM~algebra ${\cal M}(R,F)$
is a linear span of the set of
{\em Schur functions} $s_\lambda(M)$
\be
s_0(M):= 1,\qquad s_\lambda(M):=\Tr{1\dots k}
(M_{\overline 1}\dots M_{\overline k}\,
\rho_R(E_\alpha^\lambda))
\quad \lambda\vdash k,\;\; k=1,2,\dots ,
\label{def:sh-f}
\ee
where $E^\lambda_{\alpha}$ is any one of the primitive idempotents
corresponding to the partition $\lambda$ (the expression
in (\ref{def:sh-f}) does not depend
on $\alpha$).
As was shown in \cite{IOP1},
${\rm Char}(R,F)$ is an abelian algebra
with respect to the multiplication in
${\cal M}(R,F)$.

Consider
a subspace ${\rm Pow}(R,F)\subset {\rm Mat}_N({\cal M}(R,F))$ which is spanned
linearly by the elements
\ba
I\, {\rm ch}(M)\,  && \forall\; {\rm ch}(M)\in {\rm Char}(R,F)\, ,
\quad \mbox{and}
\\
M^{(x^{(k)})}:= \Tr{2\dots
k}(M_{\overline 1}\dots M_{\overline k}\, \rho_R(x^{(k)}))\,  &&
\forall\; x^{(k)}\in {\cal H}_k(q)\, , \quad
k= 1,2,\dots . \label{lam-pow}
\ea
In what follows elements of the space
${\rm Pow}(R,F)$ will be shortly called {\em the quantum matrices}.
In \cite{GPS2} it was shown that the space of the quantum matrices
carries the structure of the right ${\rm Char}(R,F)$-module
and as a ${\rm Char}(R,F)$-module
it is spanned by a series of {\em quantum matrix powers of $M$}
\be
M^{\overline 0}:= I, \quad
M^{\overline 1}:= M,\quad M^{\overline k}:=
\Tr{2\dots k}
(M_{\overline 1}\dots M_{\overline k}
\,R_{k-1} \dots R_1), \quad  k=2,3,\dots \;.
\label{k-pow}
\ee
In section 4.4 of \cite{OP2} an analogue of the matrix
multiplication was introduced for the space ${\rm Pow}(R,F)$.
It was shown there
that the {\em quantum matrix multiplication}
agrees with the right ${\rm Char}(R,F)$-module structure; it
is associative (see proposition 4.12)
and, moreover it is commutative (see propositions 4.13, 4.14).
The latter result should not be surprising as all the elements
of ${\rm Pow}(R,F)$ are descendants of the only quantum matrix
$M$.\footnote{
There should be no confusion between the quantum matrix product and
the multiplication in ${\cal M}(R,F)$. The latter one is the
product of the matrix components, while the first one
is the product of the quantum matrices.}
For our purposes in this paper it is enough knowing formulae
\be
\lb{*-product}
M^{\overline{k}} = \underbrace{M*M*\dots *M}_{\mbox{\small $k$ times}}\, ,
\qquad
(I\, {\rm ch}(M))*M^{\overline{k}}=
 M^{\overline{k}}* (I\, {\rm ch}(M))\, \quad \forall\; {\rm ch}(M)\in
{\rm Char}(R,F)\, ,
\ee
where by symbol "$*$" we denote the quantum matrix product.
We also notice that for the family of the RE~algebras the product $*$
reduces to the usual matrix product.
For the detailed description of the quantum matrix
multiplication the reader is referred to
\cite{OP2}.

The main result of our previous paper \cite{GPS2} is the
Cayley-Hamilton  theorem for the $GL(m|n)$ type QM~algebras
(see theorem I-10).
For its compact formulation and for later convenience
we introduce a shorthand notation for the following
Young diagrams (partitions)
\be
\raisebox{2mm}{$
\raisebox{5pt}{\footnotesize $r$ boxes}\!
\raisebox{5pt}{
$\left.\rule{0pt}{6mm}\right\{$}\,
\begin{tabular}{|ccc|c|}
\multicolumn{3}{c}{$\stackrel{\mbox{\footnotesize
$p$ boxes}}{\overbrace{\hspace*{16mm}}}$}
&\multicolumn{1}{c}{}\\
\hline
&\dots&&\vdots\\ \cline{4-4}
& & &\multicolumn{1}{c}{}\\ \cline{1-3}
\multicolumn{2}{|c|}{\dots}&\multicolumn{2}{c}{}\\
\cline{1-2}
\multicolumn{2}{c}{\raisebox{-3mm}{$\stackrel{
\underbrace{\hspace*{14mm}}}{\mbox{\footnotesize
$k$ boxes}}$}}&\multicolumn{2}{c}{}\\
\end{tabular}
\raisebox{12pt}{$\left.\rule{0pt}{4mm}\right\}$}
\,\raisebox{12pt}{\footnotesize $l$ boxes}$}\hspace{-3mm}
=
\left((p+1)^l,p^{(r-l)},k\right) =: [r|p]^l_k\, .
\label{lam-spec}
\ee
Here the indices $k$ and $l$ take values
$l=0,\dots ,r,$~$k=0,\dots ,p$.
If one of the indices $k$ or $l$ takes zero value, we
will omit it in the notation, e.g.,  $[r|p]^0_k = [r|p]_k$.

\begin{teor}
\lb{CH-theorem}
{\rm\bf (Cayley-Hamilton identity)~}
In the setting i)--iv) and v)--a,b)
the quantum matrix $M$ composed of the generators of the $GL(m|n)$ type
QM~algebra ${\cal M}(R,F)$ fulfills the characteristic identity
\be
\sum_{i=0}^{n+m}  M^{\overline{ m+n-i}}\,
\sum_{k=\max\{0,i-n\}}^{\min\{i,m\}} (-1)^k\,
q^{2k-i}\,
s_{[m|n]^k_{i-k}}(M)
\,\equiv 0\, .
\label{super-ch}
\ee
\end{teor}

\section*{Acknowledgment}
The authors are grateful to Alexei Davydov,
Nikolai Iorgov, Alexei Isaev,
Issai Kantor,
Hovhannes Khudaverdian, Dimitry Leites, Alexander Molev, Vladimir Rubtsov,
and Vitaly Tarasov
for numerous fruitful discussions and valuable remarks.
PP and PS gratefully acknowledge the support from RFBR grant No. 05-01-01086.
The work of PP was also
partly supported by the grant of the Heisenberg-Landau foundation.

\section{Structure of the characteristic subalgebra}
\label{sec2}

Consider the graded ring $\Lambda$ of symmetric functions in
countably many variables. A $\Bbb Z$-basis of $\Lambda$ is given
by the Schur symmetric functions $s_\lambda$, $\lambda\vdash n$, for
$n\geq 0$ (we adopt definitions and notation of ref.\cite{Mac},
sections 1.2 and 1.3).

It is not accidental that the similar notation $s_\lambda(M)$ is
assigned to the elements (\ref{def:sh-f}) of the characteristic
subalgebra of the Hecke type quantum matrix algebra
${\cal M}(R,F)$. Indeed, consider the additive map
\be
\lb{homo-map} \Lambda \ni s_\lambda\; \mapsto\;
s_{\lambda}(M)\in {\rm Char}(R,F)\subset {\cal M}(R,F)\;
\mbox{\em (Hecke type)}\, .
\ee
Our first main result is as follows.
\begin{teor}\label{teor:hom-ring}
In the setting i)--iv) and (\ref{semisimple})
the additive map (\ref{homo-map}) naturally extends to the homomorphism of rings.
\end{teor}
\addtocounter{theorem}{-1}

A proof of the theorem is given in the subsection \ref{sec2.1}.

In the subsection \ref{sec2.2} we derive some bilinear relations
for the Schur symmetric functions $s_\lambda\in \Lambda$. These
relations are necessary for the derivations in section \ref{sec3}.

\subsection{Littlewood-Richardson multiplication formula for
$s_{\lambda}(M)$}
\label{sec2.1}

We will prove the theorem \ref{teor:hom-ring} by a direct
calculation. To this end we adopt its alternative formulation

\begin{teor}{\hspace{-5pt}\bf $'$}
\label{t-1}
Let ${\cal M}(R,F)$ be a Hecke type QM~algebra
generated by the components of matrix $M$. Assume that condition
(\ref{semisimple}) on the algebra parameter $q$ is satisfied.
Then, the multiplication in the corresponding characteristic subalgebra
${\rm Char}(R,F)$ is described by the relations
\be
s_\lambda(M)s_\mu(M)  = \sum_{\nu\vdash (k+n)}
c_{\lambda\mu}^{\;\nu}s_\nu(M),
\label{lr}
\ee
where $s_{\lambda}(M), s_\mu(M) \in {\rm Char}(R,F)$
are the Schur functions (\ref{def:sh-f}),
and $c_{\lambda\mu}^{\;\nu}$ are the Littlewood-Richardson
coefficients (see, e.g., \cite{Mac}, section 1.9).
\end{teor}

\noindent{\bf Proof.~} Since the cases
$m=0$ or $k=0$ in (\ref{lr}) are trivial, we assume
$m\geq 1$ and $k\geq 1$.

Let us first prove the relation (\ref{lr})
for the case  $\mu = (1^k)$ is a single column diagram.
In that case it reads
\be
s_\lambda(M)s_{(1^k)}(M)  =
{\sum_{\nu\supset\lambda\atop\nu\vdash (k+n)}}^{\!\!\!\!\!'}
s_\nu(M).
\label{lr-1}
\ee
Here $\supset$ denotes  the inclusion relation on the set of standard
Young tableaux (see section I.2.1) and
the summation $\sum'$ is taken
only over those diagrams $\nu$
whose set theoretical difference with $\lambda$
is a vertical strip (for terminology see \cite{Mac}, section 1.1).

For single column diagrams $(1^k)$, $k=2,3,\dots,$ their corresponding
primitive idempotents $E^{(1^k)}$ satisfy the well known
iterative relations  (see, e.g. \cite{TW}, lemma 7.2, or
\cite{GPS1}, section~2.3)
\be
E^{(1)} = 1, \qquad
E^{(1^k)} = \frac{(k-1)_q}{k_q}\,
E^{(1^{k-1})}\,
\Bigl(\frac{q^{k-1}}{(k-1)_q}\, 1\, -\, \sigma_{k-1}\Bigr)
E^{(1^{k-1})} ,
\label{q-anti}
\ee
where we use notation of the section I.2.2.
We shall apply these relations for a derivation of eq.(\ref{lr-1}).
Consider the following chain of transformations
\begin{eqnarray}
\nonumber
s_\lambda(M)s_{(1^k)}(M)\, =\, \Tr{1\dots n+k}
\Bigl[\rho_R(E^\lambda_\alpha)\, \rho_R(E^{(1^k)\uparrow n})\,
M_{\overline{1}}\dots M_{\overline{n+k}}\Bigr]
\hspace{42.5mm}
\\[1mm]
\nonumber
= {(k-1)_q\over k_q}\,
\Tr{1\dots n+k}
\Bigl[\rho_R(E^\lambda_\alpha\, E^{(1^{k-1})\uparrow n})
\Bigl(\frac{q^{k-1}}{(k-1)_q}I-R_{n+k-1}
\Bigr)
\rho_R(E^{(1^{k-1})\uparrow n})\,
M_{\overline{1}}\dots M_{\overline{n+k}}\Bigr]
\hspace{1.3mm}
\\[1mm]
\nonumber
= {(k-1)_q\over k_q}\,
\Tr{1\dots n+k}
\Bigl[\rho_R(E^\lambda_\alpha\, E^{(1^{k-1})\uparrow n})
\Bigl(\frac{q^{k-1}}{(k-1)_q}I-R_{n+k-1}
\Bigr)
M_{\overline{1}}\dots M_{\overline{n+k}}\Bigr] =\dots
\hspace{14.3mm}
\\[1mm]
=
{1\over k_q}\,
\Tr{1\dots n+k}
\Bigl[\rho_R(E^\lambda_\alpha)
\Bigl(qI-R_{n+1}\Bigr)
\dots
\Bigl(\frac{q^{k-1}}{(k-1)_q}I-R_{n+k-1}\Bigr)
M_{\overline{1}}\dots M_{\overline{n+k}}\Bigr].\hspace{22mm}
\label{ss-1}
\end{eqnarray}
Here in the first line we substitute definition (\ref{def:sh-f})
for the Schur functions and use eq.(I.3.19) for $s_{(1^k)}(M)$
(the notation $E^{\mu\uparrow n}_{\beta}$ is described in lemma I.6).
We remind that this expression is independent of the choice of
index $\alpha$ labelling the primitive
idempotents $E^\lambda_\alpha\in{\cal H}_n(q)$.
In the second line we apply formula (\ref{q-anti}) (recall that
$R_i = \rho_R(\sigma_i)$).
In the third line we use relations (\ref{mk-r}) to permute
the term $\rho_R(E^{(1^{k-1})\uparrow n})$ with the
product of matrices $M$, then apply cyclic property
of the R-trace to move $\rho_R(E^{(1^{k-1})\uparrow n})$
to the leftmost position,
and take into account the commutativity of the idempotents
$E^{(1^{k-1})\uparrow n}$ and
$E^{\lambda}_\alpha$. Repeating these transformations $(k-1)$
times we eventually obtain the last line expression.

\smallskip
Let us denote the argument of the R-traces in
(\ref{ss-1}) as
\ba
Q(R)&:=& \rho_R(E^\lambda_\alpha)\, X_{n+1}\, ,
\label{Q(R)}
\\
\mbox{where}\hspace{18mm}&&
\nonumber
\\[-3mm]
X_i &:=&
\Bigl({q^{i-n}\over (i-n)_q}I-R_i\Bigr)
\Bigl(\frac{q^{i-n+1}}{(i-n+1)_q}I-R_{i+1}\Bigr)\dots
\Bigl(\frac{q^{k-1}}{(k-1)_q}I-R_{n+k-1}\Bigr).
\hspace{8mm}
\label{X}
\ea
We notice that in view of relations (\ref{mk-r}) and the cyclic
property of the R-trace one can perform
cyclic permutations of factors in $Q(R)$
without altering the expression (\ref{ss-1}).
We shall use this cyclic invariance in order to transform
$Q(R)$ to a suitable form.

The strategy of the transformation is as follows.
We use a sequence of resolutions of the idempotent
$E^\lambda_\alpha\in {\cal H}_n(q)$ ($\lambda\vdash n$)
in terms of idempotents $E^\nu_\beta\in {\cal H}_{n+i}(q)$
($\nu\vdash (n+i)$, $i\geq 1$)
described in (I.2.21)
\be
E^\lambda_\alpha =\sum_{\nu\supset\lambda\atop \nu\vdash (n+i)}
\sum_{\beta:\atop\beta\supset\alpha}
E^\nu_\beta.
\label{p-razl-1}
\ee
We successively increase $i$ in (\ref{p-razl-1}) from $2$ to $k$
and evaluate the factors $\bigl(q^{i-1}/(i-1)_q I - R_{n+i-1}\bigr)$
in $Q(R)$ on the idempotents $\rho_R(E^\nu_\beta)$
\be
\rho_R(E^\nu_\beta)\Bigl({q^{i-1}\over (i-1)_q}I - R_{n+i-1}\Bigr)\,
\stackrel{\circlearrowright}{=}
{(\ell_{n+i-1}\, +\, i \, - \, 1)_q\over
(i-1)_q\, (\ell_{n+i-1})_q }\, \rho_R(E^\nu_\beta) .
\lb{evaluate}
\ee
Here $\ell_j := c(j) - c(j+1)$ denotes the difference of the contents
of boxes with numbers $j$ and $(j+1)$ in the standard tableau
\tbl{\nu}{\beta} (for definitions see section I.2.1); the symbol
"$\stackrel{\circlearrowright}{=}$" means equality modulo cyclic
permutation of factors.

The evaluation rule can be argued as follows.
Observe that the relations
\be
E_\beta^\nu\,\sigma_{j}\equiv E_\beta^\nu
\Bigl(\sigma_{j} + \frac{q^{-\ell_{j}}}{(\ell_{j})_q}1 \Bigr)-
\frac{q^{-\ell_{j}}}{(\ell_j)_q}\,E_\beta^\nu
={(\ell_j+1)_q\over (\ell_j)_q} E^\nu_{\beta\,\pi_j(\beta)}
-\frac{q^{-\ell_{j}}}{(\ell_j)_q}
\,E_\beta^\nu , \quad 1\leq j\leq n+i-1,
\label{me-t}
\ee
are satisfied in the algebra ${\cal H}_{n+i}(q)$ (see (I.2.16)).
Here the symbol $E^\nu_{\beta\,\pi_j(\beta)}$ stands for the off-diagonal
matrix unit labelled by the pair of standard
Young tableaux \tbl{\nu}{\beta} and \tbl{\nu}{\pi_j(\beta)},
where the tableau  \tbl{\nu}{\pi_j(\beta)}
is obtained
from the tableau   \tbl{\nu}{\beta}
by the permutation $\pi_j$ of boxes $j$ and $(j+1)$.
If \tbl{\nu}{\pi_j(\beta)} is non-standard the term with
$E^\nu_{\beta\,\pi_j(\beta)}$ is absent in (\ref{me-t}).

Now, transform the expression
~$\rho_R(E^\nu_\beta) R_{n+i-1}=\rho_R(E^\nu_\beta\sigma_{n+i-1})$~
in the left hand side of (\ref{evaluate})
with the use of eq.(\ref{me-t}).
In $Q(R)$ the contribution of the off-diagonal
matrix unit $\rho_R(E^\nu_{\beta\,\pi_j(\beta)})$
vanishes by virtue of the cyclic invariance.
Indeed,
\ba
\nonumber
\rho_R(E^\nu_{\beta\,\pi_j(\beta)})\, X_{n+i}\, =\,
\rho_R(E^{\nu'}_{\beta'} E^\nu_{\beta\,\pi_j(\beta)})\,
X_{n+i}
\stackrel{\circlearrowright}{=}
\rho_R(E^{\nu}_{\beta\,\pi_j(\beta)})\, X_{n+i}\,\,
\rho_R(E^{\nu'}_{\beta'})\hspace{10mm} &&
\\
\lb{calcul}
=\,\rho_R(E^{\nu}_{\beta\,\pi_j(\beta)}E^{\nu'}_{\beta'})\, X_{n+i}
\, =\, 0 .&&
\ea
Here the idempotent $E^{\nu'}_{\beta'}$ corresponds to the standard
tableau \tbl{\nu'}{\beta'} obtained from the tableau \tbl{\nu}{\beta}
by removing the box with the number $(n+i)$. The first and the last
equalities in (\ref{calcul})  are consequences of eq.(\ref{p-razl-1})
and the  multiplication table for the matrix units (I.2.7).
In the second equality we made the cyclic permutation of terms
which is allowed in $Q(R)$. The factors $\rho_R(E^{\nu'}_{\beta'})$
and $X_{n+i}$ are built of the mutually commuting R-matrices
wherefrom the third equality in (\ref{calcul}) follows.

Eventually, collecting the coefficients
at the diagonal matrix unit $\rho_R(E^\nu_\beta)$ in $Q(R)$
results in the right hand side of eq.(\ref{evaluate}).

\smallskip
So, we begin the transformation of $Q(R)$. Setting $i=2$ in
(\ref{p-razl-1}) we come to the expression
\be
Q(R) = \sum_{\nu\supset\lambda\atop \nu\vdash (n+2)}\sum_{\beta:
\atop \beta\supset \alpha}
\rho_R(E^\nu_\beta)
\Bigl(qI-R_{n+1}\Bigr)\, X_{n+2}.
\label{sum-2}
\ee
For our calculation we have to specify an explicit way of enumeration
of the Young tableaux. For a given tableau \tbl{\lambda}{\alpha},
$\lambda\vdash n$, we take the index $\alpha := \{a_1, a_2, \dots a_n\}$
to be an ordered set of pairs of integers $a_i := \{x_i,y_i\}$, where
$x_i$ and $y_i$ are, respectively, the number of column and row
where the $i$-th box stands. Recall that the content of
the $i$-th box is $c(i) = x_i-y_i$ (see sec.I.2.1).

In the summation index $\beta$ in eq.(\ref{sum-2})
only the last two components vary. We shortly denote them as
$a$ and $b$, that is $\beta =\{\dots , a,b\}$.
For $a$ and $b$ in the summation (\ref{sum-2}) we have following
three possibilities.

\smallskip
\noindent {\em i)~} $a=\{x,y\}, b=\{x+1,y\}$.~
In this case $\ell_{n+1} = c(n+1) - c(n+2) = -1$.
Hence, due to relation (\ref{evaluate})
such tableaux do not contribute
to $Q(R)$.

\smallskip
\noindent {\em ii)~} $a=\{x,y\}, b=\{x,y+1\}$.~
In this case $\ell_{n+1} = c(n+1) - c(n+2) = 1$.
Hence, due to relation (\ref{evaluate}) the contributions of
such tableaux in (\ref{sum-2}) equal
\be
\lb{vertical2}
2_q\, \rho_R(E^\nu_{\{\dots ,a,b\}}) X_{n+2} .
\ee

\smallskip\noindent {\em iii)~ $a=\{x,y\}, b=\{\bar x,\bar y\}$, such that
$x\neq \bar x$ and $y\neq \bar y$.}~
In this case we combine  contributions coming from two
tableaux of the same shape with indices $\beta =\{\dots a,b\}$
and $\pi_{n+1}(\beta)=\{\dots b,a\}$.
Taking into account eq.(\ref{evaluate})
we get
\be
\lb{generic2}
\Bigl( \rho_R(E^{\nu}_{\{\dots ,a,b\}})\,
{(\ell_{n+1}\, +\, 1)_q\over (\ell_{n+1})_q}\, +\,
\rho_R(E^{\nu}_{\{\dots ,b,a\}})\,
{(\ell_{n+1}\, -\, 1)_q\over (\ell_{n+1})_q}
\Bigr)\, X_{n+2}
\ee
for the corresponding summands in (\ref{sum-2}).

Noticing that the term (\ref{vertical2}) fits the form (\ref{generic2})
with $\ell_{n+1}=1$ we can rewrite (\ref{sum-2}) as
\be
\lb{sum-2a}
Q(R)
\stackrel{\circlearrowright}{=}
\!\!\!\! {\sum_{{\nu \supset \lambda
\atop \nu\vdash (n+2)}\atop (a,b)}}^{\!\!\!\!\! '}
\Bigl(
\rho_R(E^{\nu}_{\{\dots ,a,b\}})\,
{(\ell_{n+1}\, +\, 1)_q\over (\ell_{n+1})_q}\, +\,
\rho_R(E^{\nu}_{\{\dots ,b,a\}})\,
{(\ell_{n+1}\, -\, 1)_q\over (\ell_{n+1})_q}
\Bigr)\, X_{n+2}\, ,
\ee
where the summation goes over different shape diagrams
$\nu\vdash (n+2)$ which are counted by
{\em unordered} pairs
$(a,b)$, $a=\{x,y\}$ and
$b=\{\bar x, \bar y\}$. There is an additional condition
$y\neq \bar y$
which means that in the diagram $\nu$
the boxes with numbers $(n+1)$
and $(n+2)$ can not appear in the same row.
It is this restriction which
the summation symbol
$\sum'$ refers to (c.f. (\ref{lr-1})).

For what follows it is suitable to
change our notation for $\ell_{n+1}$.
We substitute
$$
\ell_{n+1}= c(n+1) - c(n+2) \quad
\longrightarrow\quad
\ell_{ab} = (x-y)-(\bar x-\bar y)
$$
to manifest clearly the dependence on the summation variables $a$ and $b$.

\smallskip
We now proceed to the next step of the transformation. Substituting
(\ref{p-razl-1}) for $i=3$ into eq.(\ref{sum-2a})
and noticing ~$\ell_{ab}= - \ell_{ba}$~ we obtain
\be
\lb{sum-3}
\phantom{a}\hspace{-4mm}
Q(R)  \stackrel{\circlearrowright}{=}
\!\!\!\! {\sum_{{\tau \supset \lambda
\atop \tau\vdash (n+2)}\atop (a,b)}}^{\!\!\!\!\! '}\;\,
\sum_{\nu\vdash (n+3):\atop c=\nu\setminus \tau}
\Bigl(
\rho_R(E^{\nu}_{\{\dots ,a,b,c\}})\,
{\scriptstyle (\ell_{ab} + 1)_q\over\scriptstyle (\ell_{ab})_q}\, +\,
\rho_R(E^{\nu}_{\{\dots ,b,a,c\}})\,
{\scriptstyle (\ell_{ba} + 1)_q\over\scriptstyle (\ell_{ba})_q}
\Bigr) \Bigl({q^2\over 2_q}I-R_{n+2}\Bigr)X_{n+3} ,
\ee
where $c$ labels all possible complements of the diagram
$\tau\vdash (n+2)$
by the $(n+3)$-th
box. Applying relation (\ref{evaluate})
we  reduce this expression to the form
\be
\lb{sum-3a}
\phantom{a}\hspace{-3mm}
Q(R)  \stackrel{\circlearrowright}{=}
\!\!\!\! {\sum_{{\tau \supset \lambda
\atop \tau\vdash (n+2)}\atop (a,b)}}^{\!\!\!\!\! '}\;\,
\sum_{\nu\vdash (n+3):\atop c=\nu\setminus \tau}
\Bigl(
\rho_R(E^{\nu}_{\{\dots ,a,b,c\}})\,
{\scriptstyle (\ell_{ab} + 1)_q\over\scriptstyle (\ell_{ab})_q}
{\scriptstyle (\ell_{bc} + 2)_q\over\scriptstyle 2_q(\ell_{bc})_q}\, +\,
\rho_R(E^{\nu}_{\{\dots ,b,a,c\}})\,
{\scriptstyle (\ell_{ba} + 1)_q\over\scriptstyle (\ell_{ba})_q}
{\scriptstyle (\ell_{ac} + 2)_q\over\scriptstyle 2_q(\ell_{ac})_q}
\Bigr) X_{n+3} ,
\ee
Next, we observe that
the idempotents
$\rho_R(E^\nu_{\{\dots ,a,b,c\}})$ and
$\rho_R(E^\nu_{\{\dots ,b,a,c\}})$
in the expression
above
can be identified. Indeed, denoting~
$\sigma_i(\ell ):= (\sigma_i -q^{\ell}/{\ell}_q\, 1)$~ we have
\ba
\nonumber
\rho_R( E^\nu_{\{\dots ,b,a,c\}}) X_{n+3} \stackrel{\circlearrowright}{=}
\rho_R\Bigl(
\sigma_{n+1}(\ell_{ab}) E^\nu_{\{\dots ,b,a,c\}}
\Bigr) X_{n+3}\,\rho_R\Bigl(\sigma_{n+1}(\ell_{ab})\Bigr)^{-1}
\hspace{10mm}
\\
=
\rho_R\Bigl(E^\nu_{\{\dots ,a,b,c\}}
\sigma_{n+1}(-\ell_{ab}) (\sigma_{n+1}(\ell_{ab}))^{-1}
\Bigr) X_{n+3} \stackrel{\circlearrowright}{=}
\rho_R(E^\nu_{\{\dots ,a,b,c\}}) X_{n+3} ,
\ea
where the cyclic invariance together with relations (I.2.13), (I.2.10)
and (\ref{calcul}) were taken into account. Thus, from now on the order
of labels $a$ and $b$ makes no difference in the notation
$E^\nu_{\{\dots ,a,b,c\}}$ and we simplify it to $E^\nu_{\{\dots ,c\}}$.
Then, the expression (\ref{sum-3a}) reduces to
\be
\lb{sum-3b}
Q(R)  \stackrel{\circlearrowright}{=}
\!\!\!\! {\sum_{{\tau \supset \lambda
\atop \tau\vdash (n+2)}\atop (a,b)}}^{\!\!\!\!\! '}\;\,
\sum_{\nu\vdash (n+3):\atop c=\nu\setminus \tau}
\rho_R(E^{\nu}_{\{\dots ,c\}})\,
{ (\ell_{ac} + 1)_q\over (\ell_{ac})_q}{ (\ell_{bc} + 1)_q\over
(\ell_{bc})_q}\, X_{n+3} .
\ee
Here, noticing that $\ell_{ab}=\ell_{ac}-\ell_{bc}$,
we have transformed the coefficients at $\rho_R(E^{\nu}_{\{\dots ,c\}})$
using the q-combinatorial formula (\ref{a3}) for $k=2$ and $b_1=\ell_{ac}$,
$b_2=\ell_{bc}$ (see Appendix). The double summation is carried out with
the restriction that boxes $(n+1)$, $(n+2)$ and $(n+3)$ which are labelled
by $a$, $b$ and $c$ must be placed in different rows of the diagram
$\nu$.

Finally, we prepare the expression (\ref{sum-3b}) for the next step
calculation by collecting the summands which correspond to tableaux of the
same shape
\ba
\nonumber
Q(R)\, \stackrel{\circlearrowright}{=}
\!\!\!\!
{\sum_{{\nu\supset\lambda\atop \nu\vdash (n+3)}
\atop (a,b,c)}}^{\!\!\!\!\!\! '}
\Bigl(
\rho_R(E^{\nu}_{\{\dots ,a,b,c\}})\,
{\scriptstyle (\ell_{ac} + 1)_q\over\scriptstyle (\ell_{ac})_q}
{\scriptstyle (\ell_{bc} + 1)_q\over\scriptstyle (\ell_{bc})_q}
+\,
\rho_R(E^{\nu}_{\{\dots ,b,c,a\}})\,
{\scriptstyle (\ell_{ba} + 1)_q\over\scriptstyle (\ell_{ba})_q}
{\scriptstyle (\ell_{ca} + 1)_q\over\scriptstyle (\ell_{ca})_q}
\hspace{20mm}
\\[-2mm]
\lb{sum-3c}
+\,
\rho_R(E^{\nu}_{\{\dots ,c,a,b\}})\,
{\scriptstyle (\ell_{cb} + 1)_q\over\scriptstyle (\ell_{cb})_q}
{\scriptstyle (\ell_{ab} + 1)_q\over\scriptstyle (\ell_{ab})_q}\Bigr)
X_{n+3} .\hspace{5mm}
\ea
where the summation goes over different shape diagrams
$\nu\vdash (n+3)$ counted by {\em unordered} triples
$(a,b,c)$ such that neither pair of boxes $a$, $b$ and $c$
is placed at the same row of $\nu$.

\smallskip
Repeating the transformations described in eqs.(\ref{sum-3})--(\ref{sum-3c})
successively for $i=4,\dots ,k$ and using
q-combinatorial relations (\ref{a3}),
we eventually obtain
\be
\lb{sum-k}
Q(R)\, \stackrel{\circlearrowright}{=}
\!\!\!\!{\sum_{{\tau\supset\lambda\atop
\tau\vdash (n+k-1)}\atop (a_1,\dots,a_{k-1})}}^{\hspace{-5mm} '\hspace{4mm}}\,\,
{\sum_{\nu\vdash (n+k)\atop a_k=\nu\setminus\tau}}
\rho_R(E^{\nu}_{\{\dots ,a_k\}})
\prod_{i=1}^{k-1}
{ (\ell_{a_i a_k} + 1)_q\over (\ell_{a_i a_k})_q}\, .
\ee
Here the unordered $(k-1)$-tuples $(a_1,\dots a_{k-1})$
counting different shape diagrams $\tau\vdash (n+k-1)$
are subject to restriction that $\tau\setminus\lambda$
is a vertical strip.
The summation variable $a_k$ labels all possible complements of the diagram
$\tau\vdash (n+k-1)$
by the $(n+k)$-th
box.

Formula (\ref{sum-k}) is the $i=k$ step  analogue of the relation
(\ref{sum-3b}).
An important difference is the absence of the $X$-term
in the right hand side of the expression
(one can say that \mbox{$X_{n+k}=1$}). Therefore,
in the final expression for $Q(R)$
we have no need to distinguish between the different idempotents
$\rho_R(E^\nu_{\{\dots, a_k, \dots\}})$ ($a_k$ taking
various positions)
corresponding  to the same shape diagram $\nu\vdash (n+k)$.
Thus, the analogue of eq.(\ref{sum-3c}) reads
\be
\lb{sum-ka}
Q(R)\, \stackrel{\circlearrowright}{=}
\!\!\!\!{\sum_{{\nu\supset\lambda\atop
\nu\vdash (n+k)}\atop (a_1,\dots,a_{k})}}^{\hspace{-3.7mm} '}
\rho_R(E^{\nu}_{\{\dots\}})
\sum_{j=1}^k
\prod_{i=1\atop i\neq j}^{k}
{ (\ell_{a_i a_j} + 1)_q\over (\ell_{a_i a_j})_q}\,\, =\,\,
k_q\!\!
{\sum_{\nu\supset\lambda\atop
\nu\vdash (n+k)}}^{\!\!\!\!\! '}
\rho_R(E^{\nu}_{\{\dots\}})\, .
\ee
Here by $E^\nu_{\{\dots\}}$ an arbitrary primitive idempotent
corresponding to Young diagram $\nu$ is understood,
the summation $\sum'$ goes over all diagrams $\nu\vdash (n+k)$ such that
$\nu\setminus\lambda$ is a vertical strip,
and in the last equality we used q-combinatorial formula
(\ref{a2}) setting $\ell_{a_i a_j}= b_i-b_j$.

Substituting expression (\ref{sum-ka}) for $Q(R)$ in eq.(\ref{ss-1})
we derive formula (\ref{lr-1}), which is a particular example
of the Littlewood-Richardson
rule.

\smallskip
Now we are ready to prove the general case. To this end, let us argue
that elements $s_{(1^k)}(M)$, $k=0,1,\dots\,$, form a $\Bbb Z$-basis
of generators for the set of Schur functions. Indeed, with the help of
eqs.(\ref{lr-1}) it is easy to see that
$$
s_{(2^{k} 1^{m})}(M) = s_{(1^{(k+m)})}(M) s_{(1^{k})}(M) -
s_{(1^{(k+m+1)})}(M) s_{(1^{(k-1)})}(M) \quad \forall\; k\geq 1, m\geq 0\, .
$$
Then, using eqs.(\ref{lr-1}), elements $s_{(3^k,2^{m},1^{n})}(M)$
can be expressed as linear combinations of  monomials  of the type
$s_{(2^{l} 1^{p})}(M)s_{(1^r)}(M)$. Etc. Repeating this procedure
finitely many times one can express any Schur function $s_\lambda(M)$
as a polynomial in generators $s_{(1^k)}(M)$, $k=0,1,\dots ~$.
The explicit expressions are given by famous Jacobi-Trudi identities
(see \cite{Mac}, section 1.3).

At last, since the product of generators  $s_{(1^k)}(M)$
is described by the specification (\ref{lr-1}) of the Littlewood-Richardson
formula, the product of two arbitrary Schur functions $s_\lambda(M)$
and $s_\mu(M)$
is to be given by eq.(\ref{lr}).
\hfill\rule{6.5pt}{6.5pt} \smallskip

\subsection{Bilinear relations}
\label{sec2.2}

In this subsection we derive a series of bilinear relations for
the Schur symmetric functions $s_\lambda\in \Lambda$. By the
homomorphic map
(\ref{homo-map}) one can translate them to the characteristic
subalgebra of the Hecke type quantum matrix algebra. These relations
are used in section \ref{sec3.1} to split the characteristic
identity in the $GL(m|n)$ case into the product of two factors and,
thereby, to separate "even" and "odd" parts of the spectra of
quantum matrices.

Our derivation is based on the use of the Pl{\"u}cker relations
and we start from their short reminding (for details see \cite{Sturm}).

Consider a pair of $n\times n$ matrices $A = \|a_{ij}\|_1^n$
and $B=\|b_{ij}\|_1^n$. We denote the $i$-th row of the matrix $A$
as $a_{i*}$ and  introduce notation
\be
\lb{notat}
\det A\, :=\, [A]\, , \qquad
A\, :=\,
\left(
\begin{array}{ccccc}
a_{1*} &\dots  &a_{i*} & \dots & a_{n*}\\
1 &\dots  &i & \dots & n
\end{array}
\right),
\ee
where the latter symbol contains a detailed information on
the row content of $A$.
Namely, it says that the row $a_{i*}$ appears in the matrix $A$
at the $i$-th place
(counting downwards).

Let  us  fix a set of integer data $\{k, r_1,r_2,\dots ,r_k\}$
such that
$1\leq k\leq n$ and
$1\leq r_1<\dots < r_k\leq n$. Given these data
the Pl{\"u}cker
relation  reads
\begin{eqnarray}
[A] [B] = \sum_{1\le s_1<\dots <s_k\le n}&&
\left[
\begin{array}{ccccccccc}
a_{1*} & \dots & b_{s_1*} & \dots & b_{s_2*} &
\dots
& b_{s_k*} &\dots & a_{n*} \\
1 & \dots & r_1 & \dots & r_2 & \dots & r_k & \dots
&  n
\end{array}
\right]\times \nonumber\\
&&\left[
\begin{array}{ccccccccc}
b_{1*} & \dots & a_{r_1*} & \dots & a_{r_2*} &
\dots
& a_{r_k*} &\dots & b_{n*} \\
1 & \dots & s_1 & \dots & s_2 & \dots & s_k & \dots
&  n
\end{array}
\right] ,
\label{Pluk}
\end{eqnarray}
where the sum is taken over all possible sets $\{k,s_1,\dots ,s_k\}$.
We now apply the  Pl{\"u}cker relations
for the proof of
\begin{pred}
Let us fix four
integers~$r$,~$p$,~$l$~and~$k$,~such that
$1\le l \le r$ and  $1\le k \le p$.
Then in the ring $\Lambda$ of symmetric functions
the following bilinear relations are satisfied
(for the notation see (\ref{lam-spec}))
\be
\schf{[r|p]^l_k}\schf{[r|p]} =
\schf{[r-1|p-1]^{(l-1)}_{(k-1)}}\schf{[r+1|p+1]}+
\schf{[r|p]_k}\schf{[r|p]^l}.
\label{bil}
\ee

\end{pred}

\noindent{\bf Proof.}~
For the Schur symmetric function $s_\lambda$ corresponding to
a partition $\lambda=(\lambda_1,\dots,\lambda_p)$, the Jacobi-Trudi
relation reads
(see \cite{Mac}, section 1.3, eq.(3.4))
\be
s_\lambda = \det \|h_{\lambda_i-i+j}\|_{i,j =1}^m ,
\label{J-T}
\ee
where $m\geq p$ and the components of the matrix in the right hand
side are the complete symmetric functions (that is, the single row
Schur symmetric functions) $h_i:=s_{(i)}$.
By convention, $h_i:=0$ if $i<0$.

Substituting expressions (\ref{J-T}) into the left hand side
of relation (\ref{bil}) and using notation (\ref{notat})
we have
\begin{eqnarray}
\schf{[r|p]^l_k}\schf{[r|p]}\, =
&&
\left[
\begin{array}{cccccccc}
h_{p+1*} & h_{p*}&\dots & h_{p-l+2*}& h_{p-l*}& \dots & h_{p-r+1*} &
h_{k-r*} \\
1 & 2& \dots & l & l+1& \dots & r & r+1
\end{array}
\right]\times \nonumber\\
&&\left[
\begin{array}{cccccccc}
h_{p*} & h_{p-1*}&\dots  & h_{p-l+1*}& h_{p-l*} & \dots  & h_{p-r+1*} &
h_{-r*} \\
1 & 2& \dots & l&l+1 &\dots &  r & r+1
\end{array}
\right]\,  ,
\label{opa1}
\end{eqnarray}
where the symbol $h_{i*}:=(h_i,h_{i+1},h_{i+2}, \dots)$ is used for the
rows of the matrices appearing in the Jacobi-Trudi formula (\ref{J-T}).

Now, we transform the right hand side of eq.(\ref{opa1}) using the
Pl{\"u}cker relation
for the set of data $\{k=1,r_1=r+1\}$.
In this case most of the summands in formula
(\ref{Pluk}) vanish, since they contain
determinants of matrices with coinciding pairs of rows.
The only two contributing terms correspond to $s_1=l$ and $s_1=r+1$.
So, we get
\begin{eqnarray}
\nonumber
\schf{[r|p]^l_k}
\schf{[r|p]}&=&
\left[
\begin{array}{ccccccc}
h_{p+1*} &\dots & h_{p-l+2*}& h_{p-l*}& \dots & h_{p-r+1*} &
h_{p-l+1*} \\
1 &\dots & l & l+1& \dots & r & r+1
\end{array}
\right]\times
\\
\nonumber
&&
\left[
\begin{array}{ccccccc}
h_{p*} &\dots & h_{k-r*}& h_{p-l*}& \dots & h_{p-r+1*} &
h_{-r*} \\
1 &\dots & l & l+1& \dots & r & r+1
\end{array}
\right]\, +
\\[2mm]
\nonumber
& &\hspace{5mm}
\left[
\begin{array}{ccccccc}
h_{p+1*} &\dots & h_{p-l+2*}& h_{p-l*}& \dots & h_{p-r+1*} &
h_{-r*} \\
1 &\dots & l & l+1& \dots & r & r+1
\end{array}
\right]\times
\\
\label{opa2}
&&\hspace{5mm}
\left[
\begin{array}{ccccccc}
h_{p*} &\dots & h_{p-l+1*}& h_{p-l*}& \dots & h_{p-r+1*} &
h_{k-r*} \\
1 &\dots & l & l+1& \dots & r & r+1
\end{array}
\right]
\end{eqnarray}
which, by the Jacobi-Trudi relations,
is exactly the right hand side of the eq.(\ref{bil})
(to represent the first summand in the right hand side of (\ref{opa2})
as a product of two Schur functions one has to move $(r+1)$-th
row in its first factor up to the $(l+1)$-th place, and
$l$-th row in its second factor down to the $r$-th place).
\hfill \rule{6.5pt}{6.5pt}

\section{Various presentations of the Cayley-Hamilton
identity}
\lb{sec3}
In this section we derive three alternative expressions for the
characteristic identity (\ref{super-ch}).

In subsection  \ref{sec3.1} we use the results of section 2 to present
the characteristic identity for the $GL(m|n)$ type QM~algebra as a
product of two factors of orders $m$ and $n$. The factorization allows us
to introduce separately the sets of "even" and "odd" eigenvalues
for the quantum matrix $M$ of generators of the algebra.

In the subsection \ref{sec3.2}
we derive two other forms
of the Cayley-Hamilton identity.
They are written in terms of symmetric and skew-symmetric
powers of the quantum matrix  $M$, respectively.
The coefficients of these
identities are elements of the
characteristic subalgebra and we find their expressions in terms of
the Schur functions $s_{\lambda}(M)$, and in terms of the
eigenvalues of $M$.
For the case of supermatrices these two expressions for the
characteristic identity
were first derived in \cite{KT2,T}.

\subsection{Separation of "even" and "odd" spectral values}
\lb{sec3.1}

{}From the condition v) in the definition of the $GL(m|n)$ type
QM~algebra it follows immediately that
for all Young diagrams $\lambda$ containing the diagram
$((n+1)^{m+1})=[m+1|n+1]$
their corresponding Schur functions $s_\lambda$
belong to  the kernel of the
homomorphism (\ref{homo-map})
\be
s_\lambda\mapsto s_{\lambda}(M) = 0  \quad \forall\; \lambda:\;
((n+1)^{m+1})\subset\lambda \, .
\ee
Therefore, the image of bilinear relations
(\ref{bil}) with $r=m$, $p=n$ in the characteristic subalgebra
of the $GL(m|n)$ type QM~algebra reduces to
\be
s_{[m|n]_k^l}(M)\,s_{[m|n]}(M) =
s_{[m|n]_k}(M)\,s_{[m|n]^l}(M)\, \quad \forall\; k, l:\;
0\leq k\leq n, \; 0\leq l\leq m\, .
\label{bil-simp}
\ee
We shall use these relations to
factorize the characteristic
polynomial (\ref{super-ch}).
To this end we
multiply the identity (\ref{super-ch})
by the Schur function $\schf{[m|n]}(M)$ from the right and apply
eqs.(\ref{bil-simp}). The resulting expression reads
\be
\sum_{i=0}^{m+n}M^{\overline{m+n-i}}
\sum_{k={\max(0,i-n)}}^{\min(i,m)}(-q)^k
\schf{[m|n]^k}(M)\,\, q^{k-i}\schf{[m|n]_{(i-k)}}(M)\, \equiv\, 0\, .
\label{prom}
\ee
With the use of relations (\ref{*-product})
it can be immediately
turned into the quantum matrix product of two factors.
\begin{teor}
{\rm\bf (Cayley-Hamilton identity in a factorized form)~}
In the assumptions of theorem \ref{CH-theorem}
the identity (\ref{super-ch}) implies
\be
\Big(\sum_{k=0}^m (-q)^k\, M^{\overline{m-k}}
\schf{[m|n]^k}(M)
\Big) *\Big(\sum_{r=0}^n q^{-r}\, M^{\overline{n-r}}
\schf{[m|n]_r}(M)\Big)\equiv 0\, .
\label{factor-ch}
\ee
The identities (\ref{super-ch}) and (\ref{factor-ch})
are equivalent iff the Schur function $s_{[m|n]}(M)$
is invertible (i.e., in case if all conditions i)--v) are satisfied).
\end{teor}

The factorization suggests a natural parameterization for the
characteristic subalgebra. Namely, assuming that the conditions
i)--v) on the $GL(m|n)$ type QM~algebra ${\cal M}(R,F)$
are satisfied we consider a homomorphic map
from the characteristic subalgebra ${\rm Char}(R,F)$
into the algebra ${\Bbb C}[\mu,\nu]$ of polynomials in
two sets of (mutually commuting) variables $\mu:=\{\mu_i\}_{1\le i\le m}$
and $\nu:=\{\nu_j\}_{1\le j\le n}$. The map
${\rm Char}(R,F)\rightarrow {\Bbb C}[\mu, \nu]:
s_\lambda(M)\mapsto s_\lambda(\mu,\nu)$
called the {\em parameterization map}
is given by relations\footnote{Here we implicitly assume
the algebraic independence of the elements
$\frac{\schf{[m|n]^k}(M)}{\schf{[m|n]}(M)}$, $1\leq k\leq m$, and
$\frac{\schf{[m|n]_r}(M)}{\schf{[m|n]}(M)}$, $1\leq r\leq n$.}
\begin{eqnarray}
\textstyle\frac{\schf{[m|n]^k}(M)}{\schf{[m|n]}(M)}
&\mapsto &
{\textstyle \frac{\schf{[m|n]^k}(\mu,\nu)}{\schf{[m|n]}(\mu,\nu)}}
\, :=\!\!\!\!
\sum_{1\le i_1<\dots <i_m\le m}\!\!\!\! q^{-k}\mu_{i_1}\dots
\mu_{i_k} =e_k(q^{-1}\mu)\, , \quad 1\le k\le m\, ,
\label{def:mu}
\\[1mm]
\textstyle \frac{\schf{[m|n]_r}(M)}{\schf{[m|n]}(M)}
&\mapsto&
{\textstyle \frac{\schf{[m|n]_r}(\mu,\nu)}{\schf{[m|n]}(\mu,\nu)}}
\, :=\!\!\!\!
\sum_{1\le j_1<\dots <j_r\le n}\!\!\!\!(-q)^r\nu_{j_1}\dots\nu_{j_r}
=e_r(-q\nu)\, ,\quad 1\le r\le n\, .
\label{def:nu}
\end{eqnarray}
Here $e_k(\cdot)$ denotes the specialization of the elementary symmetric
function $e_k\in \Lambda$ to the elementary symmetric polynomial in finitely
many variables --- the arguments of $e_k(\cdot)$. The powers of the parameter
$q$ are introduced in order to get the simple form of the identity
(\ref{factor2-ch}) below.

Note that for the above parameterization we need assuming an invertibility
of the Schur function $s_{[m|n]}(M)$ (see condition v)-c)).
As we shall see in section \ref{sec4}, relations (\ref{def:mu}),
(\ref{def:nu}) define consistently the homomorphism of the  characteristic
subalgebra  ${\rm Char}(R,F)$ to a subalgebra of the supersymmetric polynomials
in variables $\{q^{-1}\mu_i\}$ and $\{-q\nu_j\}$ (see the definition in
section 4) \footnote{The characteristic subalgebra augmented by the
inverse Schur function $(s_{[m|n]}(M))^{-1}$
is parameterized by rational functions in  $\mu_i$ and $\nu_j$
which are symmetric (separately) with respect to variables in the subsets
$\mu$ and $\nu$ (see eqs.(\ref{def:mu}), (\ref{def:nu}) and proposition
\ref{proposition13}).}.
\smallskip

Now, it is straightforward to derive
a completely factorized formula for the characteristic polynomial
(\ref{super-ch}).
Namely,
the parameterization map defines naturally a left ${\rm Char}(R,F)$-module
structure on the algebra ${\Bbb C}[\mu,\nu]$. We shall use this structure to
construct {\em completion} of the space of quantum matrices:
$$
\overline{\rm Pow}(R,F):= {\rm Pow}(R,F)
\raisebox{-4pt}{$\bigotimes\atop
{\rm Char}(R,F)$} {\Bbb C}[\mu,\nu]\, .
$$
The quantum matrix product for the completed space $\overline{\rm Pow}(R,F)$
is given by formula
$$
(N\!\!\raisebox{-4pt}{$\bigotimes\atop {\rm Char}(R,F)$}\! x)*
(K\!\!\raisebox{-4pt}{$\bigotimes\atop {\rm Char}(R,F)$}\! y) :=
(N*K)\!\!\raisebox{-4pt}{$\bigotimes\atop {\rm Char}(R,F)$}\!
(x y)\, \;\;\;\forall\,
N,K\in {\rm Pow}(R,F)\, , \;
x, y \in {\Bbb C}[\mu,\nu]\, .
$$
It is associative as well as commutative (see discussion below
eq.(\ref{k-pow})).

Finally, notice that the characteristic identities (\ref{super-ch}) and
(\ref{factor-ch})
are written in the algebra ${\rm Pow}(R,F)$. When passing to the completed
algebra $\overline{\rm Pow}(R,F)$ we can apply
substitutions (\ref{def:mu}) and (\ref{def:nu}) in
the characteristic polynomial (\ref{factor-ch})
and, thus, we turn it to a completely factorized form
\be
(\schf{[m|n]}(M)\, I)^{*2} * \prod_{i=1}^m
(M - \mu_i I) *
\prod_{j=1}^n
( M - \nu_j I)
\equiv 0\, .
\label{factor2-ch}
\ee
Here all products are understood as the quantum matrix products.

The above totally factorized form of the Cayley-Hamilton theorem confirms
interpretation of the indeterminates
$\{\mu_i\}$ and $\{\nu_j\}$ as, respectively, {\em "even"
and "odd" eigenvalues of the quantum supermatrix $M$}.

\begin{remark}
{\rm
Consider the  Schur functions
~${\schf{[m|n+1]}(M)}$~ and
~${\schf{[m+1|n]}(M)}$~
standing at the unity matrix $I\equiv M^{\overline{0}}$
in the two factors of the characteristic polynomial (\ref{factor-ch}).
In the recent paper \cite{KhV}  it was shown that for the classical ($q=1$)
supermatrices
the ratio of these two Schur functions
(denoted there, respectively, as ${\rm Ber}^+(M)$ and ${\rm Ber}^-(M)$)
gives an invariant expression for the Berezinian of the supermatrix $M$
(see theorems 2 and 3 of \cite{KhV})
\be
\lb{ber}
{\rm Ber} (M) = {\schf{[m|n+1]}(M)}/\schf{[m+1|n]}(M)\, .
\ee
By the parameterization formulas (\ref{def:mu}), (\ref{def:nu}),
$$
{\rm Ber}(M) \mapsto (-1)^n q^{-(m+n)}\, {\textstyle \prod_{i=1}^m \mu_i /
\prod_{j=1}^n \nu_j}\, ,
$$
which supports an extension of the notation (\ref{ber}) to the case of the
quantum matrix $M$.
It would be desirable to relate formula (\ref{ber}) with
definitions of the quantum Berezinian suggested in \cite{N,LS}.

Next, assuming invertibility of the Schur function $\schf{[m|n]}(M)$, consider
a following combination
\be
\lb{D(M)}
\det(M)\, :=\,
{\schf{[m|n]^m_n}(M)/ \schf{[m|n]}(M)}\, =\,
{{\schf{[m|n+1]}(M)}\, \schf{[m+1|n]}(M)/ \schf{[m|n]}(M)^2}\, .
\ee
Here the second equality is due to the relation (\ref{bil-simp}).
It immediately follows from the identities (\ref{super-ch}), (\ref{factor-ch})
that the quantum matrix $M$ can be inverted provided the element $\det(M)$
is invertible. In view of the parametric formula
$$
\det(M)\mapsto (-1)^n q^{n-m}\, {\textstyle \prod_{i=1}^m \mu_i\,  \prod_{j=1}^n
\nu_j}\, ,
$$
it is natural to call $\det(M)$ a {\em determinant} of the quantum matrix $M$.
In the classical ($q=1$) supermatrix case the numerator of the first expression
in (\ref{D(M)}) ---
$\schf{[m|n]^m_n}(M)$ --- was called in \cite{KT1} a superdeterminant of the
supermatrix $M$.
}
\end{remark}

\begin{remark}
{\rm
We shall stress that the parameterization formulae
(\ref{def:mu})--(\ref{factor2-ch})
are obtained here at a formal algebraic level.
A different approach based on the representation theory
of the algebras was adopted in papers
\cite{JG,GL} (see also references therein and ref.\cite{Mudr}).
The latter approach is well applicable
for the family of RE~algebras, in which case the
characteristic subalgebra belongs to the center of the algebra
(see, e.g., \cite{I}, updated version,  section 3.2, proposition 5).
However, it seems hardly possible to apply this approach for the
QM~algebras in general.
}
\end{remark}

\subsection{Cayley-Hamilton identities for
skew-symmetric and symmetric matrix powers}
\lb{sec3.2}
We have already mentioned in the introduction that
for the Hecke type QM~algebras
the corresponding ${\rm Char}(R,F)$-module ${\rm Pow}(R,F)$
is spanned linearly by the set
$M^{\overline{k}}$, $k=0,1,2,\dots .$
The Cayley-Hamilton identity (\ref{super-ch}) then states that
${\rm Pow}(R,F)$ is {\em not a free span}
of the quantum matrix powers of $M$. In this subsection we
consider two other spanning sets for the space of quantum matrices ${\rm
Pow}(R,F)$
and derive equivalent forms of the Cayley-Hamilton identity in their
terms.

\smallskip
Consider  quantum matrices (c.f., with eq.(\ref{lam-pow}))
\ba
M^{[k|1]} &:=& \Tr{2\dots k}\Bigl(M_{\overline 1}\dots
M_{\overline k}\,
\rho_R(E^{[k|1]})\Bigr)\, ,
\\[1mm]
M^{[1|k]} &:=& \Tr{2\dots k}\Bigl(M_{\overline 1}\dots
M_{\overline k}\,
\rho_R(E^{[1|k]})\Bigr)\, ,\quad \mbox{(recall $[r|p]:=(p^r)$)}\, .
\ea
Following to A.M. Lopshits (see \cite{GGB}, p.342, or \cite{KT2,T})
we introduce
series of {\em skew-symmetric and symmetric quantum matrix powers} of the
quantum matrix $M$, respectively,
\ba
M^{\wedge 0}\, :=\, I\, ,
&&
M^{\wedge k}\, :=\, (-1)^{k-1}k_qM^{[k|1]} +(-q)^k s_{[k|1]}(M)\,I,
\quad k=1,2,\dots,
\label{ext-pow}
\\[-2mm]
\nonumber
\mbox{and}\hspace{23mm}&&
\\[-2mm]
M^{{\cal S}0}\, :=\, I\, ,
&&
M^{{\cal S}k}\, :=\,  k_q M^{[1|k]} + q^{-k}s_{[1|k]}(M)
\,I,\qquad\qquad\quad\;\,
k=1,2,\dots .
\label{sym-pow}
\ea
In \cite{IOP2} (see the
Cayley-Hamilton-Newton theorem there)
expressions for  $M^{[k|1]}$ and $M^{[1|k]}$
in terms of the quantum matrix powers of $M$
were derived.
Therefrom we calculate
\be
\label{ext-exp}
M^{\wedge k}\, =\, \sum_{r=0}^k(-q)^r M^{\overline{k-r}}
s_{[r|1]}(M)\, , \qquad
M^{{\cal S}k}\, =\, \sum_{r=0}^k q^{-r}M^{\overline{k-
r}}s_{[1|r]}(M)\, .
\ee
These relations can be  inverted with the use of the inverse
Cayley-Hamilton-Newton theorem \cite{IOP2} and the Wronski relations
(see, e.g., \cite{Mac}, eq.(2.6'))\footnote{For the elements of the
characteristic
subalgebra the Wronski relation was proved in \cite{IOP2}.}
\be
\lb{wronski}
\sum_{r=0}^k (-1)^r s_{[r|1]}\, s_{[1|k-r]}\, =\, \delta(k)\, ,
\ee
where $\delta(i):=1$ if $i=0$, and $\delta(i):=0$ otherwise.
The inverse relations read
\be
\label{ext-inv}
M^{\overline k}\,  =\,
\sum_{r=0}^k q^r M^{\wedge(k-r)}
s_{[1|r]}(M)
\, =\,
\sum_{r=0}^k (-q)^{-r}
M^{{\cal S}(k-r)} s_{[r|1]}(M)\, .
\ee
Formulae (\ref{ext-inv}) show that the space of quantum matrices
${\rm Pow}(R,F)$ is a ${\rm Char}(R,F)$-span of each one of the sets
$\{M^{\wedge k}\}_{k\geq 0}$~, $\{M^{{\cal S} k}\}_{k\geq 0}$.
We shall use them also for rewriting the
Cayley-Hamilton identity (\ref{super-ch})
in terms of the (skew-)symmetric matrix powers.
To simplify formulation
we introduce one more notation for the  Young
diagrams of a particular shape.
It is easier to explain it on the picture
$$
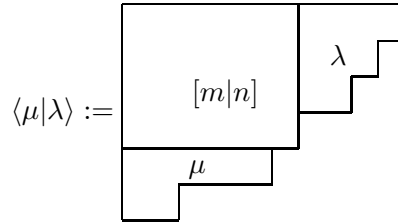
\langle\mu|\lambda\rangle :=
\begin{array}{|ccc|cccc|}\hline
 & & & & & & \\ \cline{7-7}
 & & &\multicolumn{3}{c|}{\lambda}
 &\multicolumn{1}{c}{}\\ \cline{6-6}
 &[m|n]& &\multicolumn{2}{c|}{}
 &\multicolumn{2}{c}{}\\ \cline{4-5}
 & & & \multicolumn{4}{c}{}\\ \cline{1-3}
\multicolumn{2}{|c|}{\mu}&
\multicolumn{5}{c}{}\\ \cline{2-2}
\multicolumn{1}{|c|}{\hspace*{4mm}}&
\multicolumn{6}{c}{}\\  \cline{1-1}
\end{array}
$$
that is, the Young diagram $\langle \mu|\lambda\rangle$ is a composition
of the rectangular diagram $[m|n]$ and the two diagrams
$\lambda$ and $\mu$, such that the length of $\lambda$ does not exceed $m$,
and the length of $\mu^T$ is less or equal to $n$ (we use the standard
notation from \cite{Mac}).

\begin{teor} {\rm\bf (Cayley-Hamilton identity for the
(skew-)symmetric matrix powers)}
In the assumptions of theorem \ref{CH-theorem}
the identity (\ref{super-ch})  can be written in the following
equivalent forms
\be
\sum_{k=0}^{\min\{2n, m+n\}}M^{\wedge (m+n-k)}d_k(M)\, \equiv\, 0\, ,
\qquad\mbox{or}\qquad
\sum_{k=0}^{\min\{2m, m+n\}}M^{{\cal S} (m+n-k)}f_k(M)\, \equiv\, 0\, ,
\label{ext-ch}
\ee
where we  denote
\ba
\lb{dk}
d_k(M) &:=& \textstyle{\sum_{r=\max\{0,k-n\}}^{\left[\frac{k}{2}\right]}}
(k-2r+1)_q\,s_{\langle (k-r,r)|0\rangle}(M)\, ,
\\[1mm]
\lb{fk}
f_k(M) &:=&
\textstyle{\sum_{r=\max\{0,k-m\}}^{\left[\frac{k}{2}\right]}}(-1)^{k-2r}
(k-2r+1)_q\,s_{\langle 0|(2^r,1^{k-2r})\rangle}(M)\, .
\ea
Here the symbol $\left[\frac{k}{2}\right]$ stands for the
integral part of the fraction $\frac{k}{2}$.
\end{teor}

\noindent{\bf Proof.\ } The proof of the theorem
is
a straightforward calculation on the base of
relations (\ref{ext-inv}).
We shall carry it out for the left identity in
(\ref{ext-ch}). Checking the right identity is
a similar calculation.

Substitute the expressions  (\ref{ext-inv}) for the quantum matrix powers
$M^{\overline{m+n-i}}$ in terms of the skew-symmetric powers
into the Cayley-Hamilton identity (\ref{super-ch}).
Evidently, the identity takes the form
\be
\lb{dk0}
\sum_{k=0}^{m+n} M^{\wedge{(m+n-k)}} d_k(M)\equiv 0\, ,
\ee
where the coefficients $d_k(M)\in {\rm Char}(R,F)$ are to be specified.
We shall verify  the explicit expressions (\ref{dk}) for $d_k(M)$
and refine the limits of the summation over $k$.

First of all, collecting the contributions to  $d_k(M)$
from the expressions for the matrix powers $M^{\overline{m+n-i}}$,
$0\le i\le k$, we have
$$
d_k(M) = \sum_{i=0}^kq^{k-i}s_{(k-i)}(M)
\sum_{j=\max\{0,i-n\}}
^{\min\{i,m\}}(-1)^j\, q^{2j-i}\,
s_{[m|n]_{i-j}^j}(M)\, .
$$
Then, introducing a new summation variable $r=i-j$  and changing
the order of summation over $i$ and $r$ we get
\be
d_k(M) = \sum_{r=0}^{\min\{k,n\}}(-1)^rq^{k-2r}
\sum_{i=r}^{\min
\{k,r+m\}}(-1)^i\,s_{(k-i)}(M)
s_{[m|n]_r^{i-r}}(M)\, .
\label{dk1}
\ee
Let us separately calculate the second sum in the expression above.
\begin{lem}
\label{lem:4}
For any fixed pair of integers $m$ and $n$, and for all
integers $r$ and $k$ satisfying conditions
~$0\le r\le n$,~ $r\le k\le m+n$,~ the following
equalities
\be
\sum_{i=r}^{\min\{k,r+m\}}(-1)^is_{(k-i)}
s_{[m|n]_r^{i-r}} = \left\{
\begin{array}{cl}
0\, , &k\ge  n+r+1\, ,
\\ & \\
(-1)^r
\sum_{{i= \max\{0,k-n\}}}^{\;\;\;\;\,\min\{r,k-r\}}
s_{\langle (k-i,i)|0\rangle}\, , & k\le n+r\, ,
\end{array}
\right.
\label{sum}
\ee
take place in the ring  $\Lambda$ of the symmetric functions.
\end{lem}

\noindent
{\bf Proof.\  }
Denote $\omega_{k,r}$ the expression in the left
hand side of eq.(\ref{sum}).

Consider the case $k\le r+m$.
Introducing a new summation
variable $j=k-i$ and denoting $p:=k-r$,~ $0\le p\le m$, we
rewrite the sum $\omega_{k,r}$ as
$$
(-1)^k\omega_{k,r} = \sum_{j=0}^p(-1)^js_{(j)}
s_{[m|n]_r^{p-j}} =
(-1)^ps_{(p)}s_{[m|n]_r} + \sum_{j=0}^{p-1}
(-1)^j s_{(j)}s_{[m|n]_r^{p-j}}.
$$
Applying the Littlewood-Richardson rule to the
products $s_{(j)}s_{[m|n]_r^{p-j}}$
we can gather terms in the latter expression
into two separate sums
\begin{eqnarray*}
(-1)^k\omega_{k,r} = & & \hspace*{-6mm}
(-1)^p s_{(p)}s_{[m|n]_r}+
\sum_{j=0}^{p-1}(-1)^j\sum_{t=0}^{\min\{j,n\}}
\sum_{i= \max\{0, r+t-n\}}^{\min\{r,t\}}
s_{\langle (r+t-i,\,i)|(j-t+1,\,1^{p-j-1})\rangle}\\
&&+\sum_{j=1}^{p-1}(-1)^j\sum_{t=0}^{\min\{j-1,n\}}
\sum_{i= \max\{0, r+t-n\}}^{\min\{r,t\}}
s_{\langle (r+t-i,\,i)|(j-t,\,1^{p-j})\rangle}\, .
\end{eqnarray*}
As can be easily checked,
the two triple sums in the expression above cancel each other
except for the term
$j=p-1$ in the first sum. So, we obtain
\be
(-1)^k \omega_{k,r} = (-1)^p s_{(p)}s_{[m|n]_r} +
(-1)^{p-1}\sum_{t=0}^{\min\{p-1,n\}}
\sum_{i= \max\{0, r+t-n\}}^{\min\{r,t\}}
s_{\langle (r+t-i,i)|{(p-t)}\rangle}.
\label{om}
\ee
Consider now expansion of the product
$s_{(p)}s_{[m|n]_r}$
into the sum of Schur  symmetric functions
\be
(-1)^ps_{(p)}s_{[m|n]_r} = (-1)^p
\sum_{t=0}^{\min\{p,n\}}
\sum_{i = \max\{0,\, t+r-n\}}^{\min\{r,t\}}
s_{\langle (r+t-i,\,i)|(p-t)\rangle}\, .
\label{sps}
\ee
Comparing the double sums in the right hand sides
of eqs.(\ref{om}) and (\ref{sps})
we observe that they are exactly opposite in the sign in case
$k\ge r+n+1 \Leftrightarrow p\ge n+1$, and they differ by the
term with $t=p$ in case $k\le n+r \Leftrightarrow p\le n$.
Therefore, substitution of the expression (\ref{sps}) into eq.(\ref{om})
results in formula (\ref{sum}).

The case $k\ge r+m+1$ is
treated in complete analogy with the above consideration.
\hfill\rule{6.5pt}{6.5pt}
\smallskip

Return to the proof of the theorem. By the homomorphism
(\ref{homo-map}) the statement of lemma \ref{lem:4} translates to the ring of
Schur functions $s_\lambda(M)$.
So, formula (\ref{dk1}) for $d_k(M)$  can be equivalently written as
\be
\lb{dk2}
d_k(M)=0, \;\;\;\mbox{if $k>2n$}; \qquad
d_k(M) = \!\!\!\!\sum_{r=\max\{0,k-n\}}^{\min\{k,n\}}
\!\!\! q^{k-2r}\!\!
\sum_{i=\max\{0,k-n\}}^{\min\{r,k-r\}}
s_{\langle (k-i,\,i)|0 \rangle},\;\;\;
\mbox{for $0\le k\le 2n$} .
\ee
The latter expression can be further simplified.
In case $0\le k\le n$ we have
\be
\lb{dk3}
d_k(M) =\sum_{r=0}^k q^{k-2r}\!\!
\sum_{i=0}^{\min\{r,k-r\}}\!\!\!\!
s_{\langle (k-i,\,i)|0\rangle}
=\sum_{i=0}^{[k/2]}s_{\langle (k-i,\,i)|0\rangle}
\sum_{r=i}^{k-i}q^{k-2r} =
\sum_{i=0}^{[k/2]}(k-2i+1)_q\,
s_{\langle  (k-i,\,i)|0\rangle}\, ,
\ee
where in the second equality we changed the order of summation.
In case $n<k\le 2n$ the similar calculation gives
\be
\lb{dk4}
d_k(M) =\sum_{r=k-n}^n q^{k-2r}\sum_{i=k-n}^{\min\{r,k-r\}}
s_{\langle (k-i,\,i)|0\rangle} =
\sum_{i=k-n}^{[k/2]}(k-2i+1)_q\, s_{\langle (k-i,\,i)|0 \rangle}.
\ee
Combining together the results  (\ref{dk0}), (\ref{dk2}), (\ref{dk3}) and
(\ref{dk4})
we get formulae (\ref{ext-ch}), (\ref{dk}).
\phantom{aaaa}\hfill\rule{6.5pt}{6.5pt}
\smallskip

Assuming additionally the Schur function $s_{[m|n]}(M)= d_0(M) = f_0(M)$
to be invertible we will now express the ratios
$d_k(M)/d_0(M)$ and $f_k(M)/f_0(M)$
in terms of the eigenvalues
of the quantum supermatrix $M$.
\begin{pred}
Let ${\cal M}(R,F)$ be a QM~algebra of
the $GL(m|n)$ type, that is the algebra defined by the set of
conditions i)--v) (see introduction). Then,
under the parameterization map (\ref{def:mu}), (\ref{def:nu})
we have
\begin{eqnarray}
\frac{d_k(M)}{d_0(M)} &\mapsto&
(-1)^k \sum_{r=\max\{0,k-n\}}^
{\min\{k,n\}} q^{2r}e_r(\nu)\,e_{k-r}(\nu)\, ,
\label{dkd0}
\\[1mm]
\frac{f_k(M)}{f_0(M)} &\mapsto&
(-1)^k \sum_{r=\max\{0,k-m\}}^
{\min\{k,m\}} (-q)^{-2r}e_r(\mu)\,e_{k-r}(\mu)\, .
\label{fkf0}
\end{eqnarray}
\end{pred}

\noindent
{\bf Proof.\ } We shall prove the equality (\ref{dkd0}).
The relation (\ref{fkf0}) can be checked in a similar way.

Multiplying eq.(\ref{dk1}) by $d_0(M)$ we obtain
\begin{eqnarray*}
d_k(M)d_0(M) &=&
\sum_{l=0}^{\min\{k,n\}}(-1)^lq^{k-2l}
\sum_{j=l}^{\min\{k,l+m\}}(-1)^j\,s_{(k-j)}(M)\,
\schf{[m|n]^{(j-l)}_l}(M)\,\schf{[m|n]}(M)\\
& = & \sum_{l=0}^{\min\{k,n\}}q^{k-2l}\schf{[m|n]_l}(M)
\sum_{j=0}^{\min\{k-l,m\}}(-1)^j\,s_{(k-l-j)}(M)\,
\schf{[m|n]^j}(M),
\end{eqnarray*}
where in passing to the second line we apply the bilinear
relations (\ref{bil-simp}) and shift the summation index $j\rightarrow j-l$.
The last sum in the second line
can be calculated with the use of relation (\ref{sum})
(take there $r=0$ and substitute $k\rightarrow k-l$).
The result is
\be
\lb{dkd0-1}
d_k(M)d_0(M) = \sum_{l=\max\{0,k-n\}}^{\min\{k,n\}}
q^{k-2l}\schf{[m|n]_l}(M)\,\schf{[m|n]_{(k-l)}}(M).
\ee
The parameterization formula (\ref{dkd0}) follows immediately from
the relations (\ref{def:nu}) and (\ref{dkd0-1}).
\hfill\rule{6.5pt}{6.5pt}.

\section{Spectral parameterization of the characteristic
subalgebra}
\lb{sec4}

In this section we complete the parameterization of the
characteristic subalgebra in terms  of the eigenvalues of quantum supermatrix
$M$. To this end, in the subsection \ref{sec4.1} we derive parametric
expressions for the generators $s_{(1^k)}(M)=s_{[k|1]}(M)$ and
$s_{(k)}(M)=s_{[1|k]}(M)$
and prove that the characteristic subalgebra is parameterized by the
supersymmetric polynomials.
This, in principle, solves the parameterization problem.

In the last subsection
\ref{sec4.2} we derive parameterization formula (\ref{c0}) for the Schur
function
$\schf{[m|n]}(M)$.
The latter result allows translating the condition of invertibility of
$s_{[m|n]}(M)$ into conditions on the spectral variables
$\{\mu_i\}_{1\leq i\leq m}$ and $\{\nu_j\}_{1\leq j\leq n}$.
We shall prove formula (\ref{c0}) using yet another series of bilinear relations
in the ring $\Lambda$ of symmetric functions (see lemma \ref{lem:12}).
Note that relation (\ref{c0}) is a particular case of the factorization
formula known in the theory of the supersymmetric polynomials
\cite{BR,PrT}.

\subsection{Parameterization of the single column
and the single row Schur functions}
\lb{sec4.1}

\begin{pred}\label{pred:par}
Let ${\cal M}(R,F)$ be the $GL(m|n)$ type QM~algebra
satisfying the conditions i)--v) (see introduction). Then,
the parameterization map (\ref{def:mu}), (\ref{def:nu})
assigns the following expressions to the generators
$\{s_{[k|1]}(M)\}_{k\geq 0}$
and $\{s_{[1|k]}(M)\}_{k\geq 0}$
of the characteristic subalgebra ${\rm Char}(R,F)$
\begin{eqnarray}
s_{[k|1]}(M)&\mapsto &
s_{[k|1]}(\mu,\nu)\, =\, \sum_{r=0}^ke_r(q^{-1}\mu)\,h_{k-r}(-q\nu)\, ,
\label{sk-a}\\
s_{[1|k]}(M)&\mapsto &s_{[1|k]}(\mu,\nu)
\, =\,  \sum_{r=0}^ke_r(-q\nu)\,h_{k-r}(q^{-1}\mu)\, .
\label{sk-s}
\end{eqnarray}
$$
\!\!\!\!\!\!\!\!\!\!
\mbox{Here}\quad
e_r(q^{-1}\mu):=\!\!\!\!\sum_{1\le i_1<\dots<i_r\le m}\!\!\!\! q^{-r}
\mu_{i_1}\mu_{i_2}\dots \mu_{i_r}\quad \mbox{and}\quad
h_r (-q \nu) :=\!\!\!\!\sum_{1\le i_1\le \dots\le i_r\le n}
\!\!\!\!(-q)^r \nu_{i_1}\nu_{i_2}\dots \nu_{i_r}
$$
are the elementary symmetric and complete symmetric polynomials
in $m$ and $n$ variables, respectively (\cite{Mac}, section 1.2).
\end{pred}

\noindent
{\bf Proof.\ \ }
We apply induction on $k$.
By the Littlewood-Richardson rule (\ref{lr})
$$
\schf{[m|n]}(M)\,s_{(1)}(M) = \schf{[m|n]^1}(M)+
\schf{[m|n]_1}(M).
$$
Dividing both sides of this equality by $\schf{[m|n]}(M)$ and
using relations (\ref{def:mu}),  (\ref{def:nu})
we get the parameterization formula for $s_{(1)}(M)$
$$
s_{(1)}(\mu,\nu) = s_{[1|1]}(\mu,\nu)
= e_1(q^{-1}\mu) +e_1(-q\nu)\, ,
$$
which can be equivalently written as
$$
s_{[1|1]}(\mu,\nu) =  e_1(q^{-1}\mu) + h_1(-q\nu)\, , \quad\mbox{or as}
\quad s_{[1|1]}(\mu,\nu) =  e_1(-q\nu) + h_1(q^{-1}\mu)\, .
$$
These formulae are nothing but the eqs.
(\ref{sk-a}) and (\ref{sk-s}) in case $k=1$.

Now, assuming the relations (\ref{sk-a}) and
(\ref{sk-s}) are valid for all values of the index  $1\leq k<p$
we shall prove them
for $k=p$. For definiteness, we  check
the  eq.(\ref{sk-s}). The eq.(\ref{sk-s}) is worked out
similarly.

Let us write down the  image
of  the relation (\ref{sum}) in the characteristic subalgebra,
the case $r=0$, $k=p$ :
\be
\sum_{i=0}^{\min\{p,m\}} (-1)^i s_{(p-i)}(M)\, s_{[m|n]^i}(M)\,
=\, \theta(n-p)\, s_{[m|n]_p}(M)\, .
\label{sum-spec}
\ee
Here $\theta(i):= 0$ if $i<0$, and  $\theta(i):= 1$ otherwise.
Substituting the parametric expressions (\ref{def:mu}) and (\ref{def:nu})
for $s_{[m|n]^i}(M)/s_{[m|n]}(M)$ and $s_{[m|n]_p}(M)/s_{[m|n]}(M)$
into (\ref{sum-spec}) we find
\be
\lb{kusok}
s_{[1|p]}(\mu,\nu) = e_{p}(-q\nu) -
\sum_{i=1}^{p}(-1)^{i}s_{[1|p-i]}(\mu, \nu)\,e_{i}(q^{-1}\mu).
\ee
Now,  using the induction assumption we substitute expressions~ (\ref{sk-s})~
for~ the~ elements $s_{[1|p-i]}(\mu,\nu)$,~ $1\leq i\le p$,~
into (\ref{kusok}) and calculate
\ba
\nonumber
s_{[1|p]}(\mu,\nu) &=& e_{p}(-q\nu) -\sum_{j=0}^{p-1}
e_j(-q\nu)\sum_{i=1}^{p-j}
\left( (-1)^i h_{p-j-i}(q^{-1}\mu)\,e_{i}(q^{-1}\mu)\right)
\\[1mm]
\lb{opsa}
&=&
e_{p}(-q\nu) + \sum_{j=0}^{p-1} e_j(-q\nu)\,
h_{p-j}(q^{-1}\mu) =  \sum_{j=0}^{p}e_j(-q\nu)\,
h_{p-j}(q^{-1}\mu)\, ,
\ea
where in passing to the second line we
used the Wronski relations (\ref{wronski})
for the substitution
$$\sum_{i=1}^{p-j}
(-1)^i h_{p-j-i}(q^{-1}\mu)\,e_{i}(q^{-1}\mu)\, =\, -h_{p-j}(q^{-1}\mu)\, .
$$
Calculation (\ref{opsa})
completes the inductive proof of the eq.(\ref{sk-s}).
\hfill
\rule{6.5pt}{6.5pt}
\smallskip

Let us recall the definition of the supersymmetric polynomials
(see, e.g.,  \cite{Stem}).

\begin{opred}
Let $x=\{x_i\}_{1\le i\le m}$ and
$y=\{y_j\}_{1\le j\le n}$ be two sets of independent commutative
variables. A polynomial $p\in {\Bbb C}[x,y]$ is said to be
{\it supersymmetric} if
\begin{itemize}
\item[{\it a)}] $p$ is invariant under permutations of $x_1,\dots ,x_m$;
\item[{\it b)}] $p$ is invariant under permutations of $y_1,\dots ,y_n$;
\item[{\it c)}] upon substituting $x_1 = y_1 = t$
in $p$, the resulting polynomial
does not depend on $t$.
\end{itemize}
An algebra of the supersymmetric polynomials is further denoted as
$T[x,y]$.
\end{opred}

Obviously, the polynomials $s_{[k|1]}(\mu,\nu)$ and $s_{[1|k]}(\mu,\nu)$
given by eqs.(\ref{sk-a}), (\ref{sk-s}) satisfy the conditions {\it a)}
and {\it b)} of the above definition
with respect to  variables $x_i=q^{-1}\mu_i$ and $y_j=-q\nu_j$.
Validity of the property {\it c)} for them results from the following statement.
\begin{lem}
Denote
$\{\mu'\}:=\{\mu\}\setminus \{\mu_1\}=\{\mu_i\}_{2\le i\le m}$~,
$\{\nu'\}:=\{\nu\}\setminus\{\nu_1\}=\{\nu_i\}_{2\le j\le n}$.
Then  the polynomials $s_{[1|k]}(\mu,\nu)$ and $s_{[k|1]}(\mu,\nu)$
satisfy expansions
\begin{eqnarray}
s_{[1|k]}(\mu,\nu) &=& s_{[1|k]}(\mu',\nu') + (q^{-1}\mu_1 - q\nu_1)
\sum_{r=0}^{k-1}\left(q^{-1}\mu_1\right)^{k-r-1}s_{[1|r]}(\mu',\nu')\, ,
\label{ss-dec}\\
s_{[k|1]}(\mu,\nu) &=& s_{[k|1]}(\mu',\nu') + (q^{-1}\mu_1 - q\nu_1)
\sum_{r=0}^{k-1}(-q\nu_1)^{k-r-1}\, s_{[r|1]}(\mu',\nu').
\label{sa-dec}
\end{eqnarray}
\end{lem}

\noindent
{\bf Proof.\ } For the elementary and complete symmetric functions
one has
$$
e_k(\mu)\, =\, e_k(\mu')+ \mu_1\,
e_{k-1}(\mu')\, , \qquad
h_k(\mu)\, =\, \sum_{r=0}^k(\mu_1)^r\,
h_{k-r}(\mu')\, .
$$
Substituting these formulae into eqs.(\ref{sk-a}) and (\ref{sk-s})
it is easy to derive formulae (\ref{ss-dec}), (\ref{sa-dec}).
\hfill
\rule{6.5pt}{6.5pt}
\smallskip

We have checked that the polynomials $s_{[k|1]}(\mu,\nu)$ and
$s_{[1|k]}(\mu,\nu)$ are supersymmetric.  Moreover, as was proved
in \cite{PrT} (see theorem (3.1) and proposition (2.3) there)
the algebra of supersymmetric polynomials $T[q^{-1}\mu,-q\nu]$
can be generated by any one of the sets
$\{s_{[1|k]}(\mu,\nu)\}_{k\geq 0}$ or
$\{s_{[k|1]}(\mu,\nu)\}_{k\geq 0}$.
Therefore, as a direct consequence of  the proposition \ref{pred:par}
we get
\begin{cor}
\label{cor:8} In the conditions of the proposition \ref{pred:par}
an image of the characteristic subalgebra ${\rm Char}(R,F)$
under the parameterization map (\ref{def:mu}), (\ref{def:nu})
is the algebra $T[q^{-1}\mu,-q\nu]$ of the supersymmetric polynomials
in the variables $\{q^{-1}\mu_i\}_{1\le i\le m}$ and
$\{-q\nu_j\}_{1\le j\le n}$.
\end{cor}

\subsection{Parameterization of the Schur function $s_{[m|n]}(M)$}
\label{sec4.2}

In a derivation of the parameterization     formula for
$s_{[m|n]}(M)$ we will use A. Kirillov's bilinear relations on the
Schur functions (see \cite{Kir,KR,Kl})
\be
\lb{kir}
s_{[m|n]}\,
s_{[m|n]}\, =\, s_{[m+1|n]}\, s_{[m-1|n]}\,  +\, s_{[m|n+1]}\,
s_{[m|n-1]}\,  \qquad\forall \; m,n=1,2,\dots \, .
\ee
To keep our presentation self-contained let us briefly describe how one can
prove them
using the Pl\"{u}cker relations and the Jacobi-Trudi
formulae. Actually,
one can derive a more extensive set of relations.
\begin{lem}
\lb{lem:12}
For any integers $a$, $b$, $m$, $n$ : $1\leq a\leq m$,~ $1\leq b\leq n$,~
the equalities
\begin{eqnarray}
s_{[a|b]}\, s_{[m|n]} & = &\!\!\!\!\!\!\sum_{k=\max\{1,a+b-n\}}^a
\!\!\!\!(-1)^{a-k}
s_{[m|n]_{a+b-k}}\, s_{[a-1|b-1]^{k-1}}
\nonumber
\\
&&+
\!\! \sum_{k=\max\{1, a+b-m\}}^b
\!\!\!\!(-1)^{b-k}
s_{[m|n]^{a+b-k}}\, s_{[a-1|b-1]_{k-1}}\hspace{10mm}
\label{sab-exp}
\end{eqnarray}
are satisfied in the ring of symmetric functions $\Lambda$.

Formulae (\ref{kir}) correspond to the choice $a=m$, $b=n$ in
eqs.(\ref{sab-exp}).
\end{lem}

\noindent
{\bf Proof.\  }
Applying the Jacobi-Trudi relation (\ref{J-T})
and using elementary properties of determinants
we can write down determinantal presentations for the Schur
functions  $s_{[m|n]}$ and $s_{[a|b]}$
\ba
\lb{odin}
s_{[m|n]}& =&
\left[
\begin{array}{ccccc}
h_{n*} & h_{n-1*}&\dots & h_{n-m+1*}& \delta_{m+1 *} \\
1 & 2& \dots & m & m+1
\end{array}
\right]\, ,
\\
\lb{dva}
s_{[a|b]} & =&
\left[
\begin{array}{cccccccc}
\delta_{1*} & \delta_{2*}&\dots & \delta_{m+1-a*}&
h_{b-m+a-1*}&\dots&h_{b-m+1*}&h_{b-m*} \\
1 & 2& \dots & m+1-a & m+2-a&\dots &m&m+1
\end{array}
\right]\, .
\ea
Here we use the matrix notation introduced in (\ref{notat}) and
the symbols $h_{i*}$ and $\delta_{i*}$ denote the following matrix
rows
$$
h_{i*}:= \left(h_i,h_{i+1},h_{i+2}, \dots\right)\, ,\qquad
\delta_{i*}:=\bigl(
0,  \dots,0, 1\hspace{-5pt}\raisebox{12pt}{$\downarrow\!\!$
\lefteqn{\raisebox{3pt}{\footnotesize\em $i$-th place}}}, 0, \dots \bigr)\, .
$$
Relations (\ref{sab-exp}) result from an application of the Pl\"{u}cker
relations
(\ref{Pluk}) for the set of data $\{k=1,r_1=m+1\}$ to the product of
determinants (\ref{odin}) and (\ref{dva}).
\hfill
\rule{6.5pt}{6.5pt}\smallskip

Now we are ready to prove the main result of this subsection.

\begin{pred}
\lb{proposition13}
Let ${\cal M}(R,F)$ be the $GL(m|n)$ type QM algebra satisfying the
conditions i)--v) (see introduction).
Then, the image of the
Schur function $\schf{[m|n]}(M)$
under the parameterization map (\ref{def:mu}), (\ref{def:nu})
is given by formula
\be
\schf{[m|n]}(M)\, \mapsto\, s_{[m|n]}(\mu,\nu) = \prod_{i=1}^m\prod_{j=1}^n
\left(q^{-1}\mu_i - q\nu_j\right).
\label{c0}
\ee
Therefore, the invertibility of the Schur function $s_{[m|n]}(M)$
implies invertibility of all factors $\left(q^{-1}\mu_i - q\nu_j\right)$
in the product (\ref{c0}) for
$1\leq i\leq m$~ and $1\leq j\leq n$.
\end{pred}

\noindent
{\bf Proof.\  }
Let us first multiply the image of the relation (\ref{kir}) in the
characteristic subalgebra ${\rm Char}(R,F)$ by $(s_{[m|n]}(M))^{-1}$
and then apply the parameterization map.
By virtue of the relations (\ref{def:mu}), (\ref{def:nu})
the resulting formula reads
\be
\schf{[m|n]}(\mu,\nu) = e_n(-q\nu)\,\schf{[m-1|n]}(\mu,\nu)
\, +\, e_m(q^{-1}\mu)\,\schf{[m|n-1]}(\mu,\nu).
\ee
Noticing that
$$
e_m(q^{-1}\mu)|_{\mu_i = 0} = 0\,
\quad \forall i=1,\dots ,m, \qquad
e_n(-q\nu)|_{\nu_j = 0} = 0\, \quad
\forall j=1,\dots ,n,
$$
we obtain for the supersymmetric polynomial $s_{[m|n]}(\mu,\nu)$
\be
\schf{[m|n]}(\mu,\nu)|_{q^{-1}\mu_i = q\nu_j} =
s_{[m|n]}(\mu,\nu)|_{\mu_i=\nu_j=0} = 0\, \quad \forall i, j:\;
1\le i\le m, \; 1\le j\le n.
\label{snm-0}
\ee
As immediately follows from the Jacobi-Trudi relation (\ref{J-T}), the Schur
function $\schf{[m|n]}(\mu,\nu)$ is a homogeneous polynomial
in the variables $\{q^{-1}\mu_i\}_{1\le i\le m}$ and $\{q\nu_j\}_{1\le j\le n}$
of the order $(m+n)$. Together with eq.(\ref{snm-0}) this implies
\be
\schf{[m|n]}(\mu,\nu) =\alpha\, \prod_{i=1}^m\prod_{j=1}^n
\left(q^{-1}\mu_i - q\nu_j\right)\, ,
\label{anz}
\ee
where $\alpha$ is a numeric factor.
To define $\alpha$, observe the following consequence of
(\ref{sk-s})
$$
s_{(k)}(\mu,\nu)|_{\mu_1 = \dots=\mu_m = 0} = e_k(-q\nu).
$$
Therefore
$$
\schf{[m|n]}(\mu,\nu)_{|_{\mu_1 = \dots=\mu_m = 0}} =
\det\Bigl(e_{n-i+j}(-q\nu)\Bigr)_{i,j=1}^m
= \Bigl(e_n(-q\nu)\Bigr)^m = \Bigl(\prod_{i=1}^n(-q\nu_i)\Bigr)^m.
$$
Comparing this result with eq.(\ref{anz})  at the point $\mu_1=\dots =\mu_m=0$,
we
find $\alpha=1$ thereby ending the proof.
\hfill\rule{6.5pt}{6.5pt}

\section*{Appendix}
\label{appen}
\addtocontents{toc}{\contentsline {section}{\numberline {} Appendix}{\pageref{appen}}}

\renewcommand{\theequation}{{\small A}.\arabic{equation}}
\setcounter{equation}0

\noindent
\rm
Here we derive  the q-combinatorial relations which are used in the proof
of theorem \ref{t-1}.

For an arbitrary set of pairwise different nonvanishing integers
$b_i$, $i=1,2,\dots ,k$,
we shall prove following relations
\ba
\lb{a1}
q^k\ -\ \prod_{i=1}^k {(b_i+1)_q\over (b_i)_q} &=&
 -\, \sum_{j=1}^k {q^{-b_j}\over (b_j)_q}
\prod_{\scriptstyle i=1\atop\scriptstyle i\neq j}^{k}
{(b_i-b_j+1)_q\over (b_i-b_j)_q}\, ,
\\[2mm]
\lb{a2}
k_q &=&
\sum_{j=1}^k\prod_{\scriptstyle i=
1\atop\scriptstyle  i\neq j}^{ k}
{(b_i-b_j+1)_q\over (b_i-b_j)_q} ,
\\[2mm]
\lb{a3}
\prod_{i=1}^k {(b_i+1)_q\over (b_i)_q} &=&
\sum_{j=1}^k {(b_j+k)_q\over k_q (b_j)_q}
\prod_{\scriptstyle i=1\atop\scriptstyle i\neq j}^{k}
{(b_i-b_j+1)_q\over (b_i-b_j)_q}
\ea
A proof is by induction on $k$. Checking  the case $k=1$
in relations (\ref{a1})--(\ref{a3}) is an easy exercise.
Now, assuming  relations (\ref{a1}) are valid for all $k\leq m$
let us transform the expression in the
left hand side of eq.(\ref{a1}) for
$k=m+1$
\ba
\nonumber
\hspace{-3mm}
q^{m+1}-\prod_{i=1}^{m+1}{(b_i+1)_q\over (b_i)_q}
&=&
\left(q^{m+1} -
q^m{(b_{m+1}+1)_q\over (b_{m+1})_q}\right)\ +\
\left(q^m - \prod_{i=1}^{m}{(b_i+1)_q\over (b_i)_q}
\right)
{(b_{m+1}+1)_q\over (b_{m+1})_q}
\\[2mm]
\lb{a4}
&=& -\
{q^{m-b_{m+1}}\over (b_{m+1})_q}\ -\
\left(
\sum_{j=1}^m {q^{-b_j}\over (b_j)_q}
\prod_{\scriptstyle i=1\atop\scriptstyle i\neq j}^{m}
{(b_i-b_j+1)_q\over (b_i-b_j)_q}\right)
{(b_{m+1}+1)_q\over (b_{m+1})_q}\, .
\ea
Here we used formula  (\ref{a1}), case
$k=m$,
for the transformation of the last term in the first line.
For further transformation we use the formula
\be
\lb{a5}
{(b_{m+1}+1)_q\over (b_{m+1})_q} = {(b_{m+1}-b+1)
_q\over (b_{m+1}-
b)_q} -
{(b)_q\over (b_{m+1})_q(b-b_{m+1})_q}\, .
\ee
Substituting $b_j$, $j=1,2,\dots m$,
for $b$  in eq.(\ref{a5}) we continue the calculation
\ba
\nonumber
\mbox{(\ref{a4})}
\hspace{-2mm}
&=&
\hspace{-2mm}
-\sum_{j=1}^m {q^{-b_j}\over (b_j)_q}
\prod_{\scriptstyle i=1\atop\scriptstyle i\neq j}^{m+1}
{(b_i-b_j+1)_q\over (b_i-b_j)_q}
\ -\
{q^{-b_{m+1}}\over (b_{m+1})_q}
\left( q^m +
\sum_{j=1}^m {q^{-b_j+b_{m+1}}\over (b_j-b_{m+1})_q}
\prod_{\scriptstyle i=1\atop\scriptstyle i\neq j}^{m}
{(b_i-b_j+1)_q\over (b_i-b_j)_q}\right) ,
\ea
and then, applying relation (\ref{a1}) with the
shifted set of integers
$b_i \rightarrow (b_i-b_{m+1})$, $i=1,2,\dots ,m$,
for the transformation of the last term we obtain
\ba
\nonumber
&=&
- \sum_{j=1}^m {q^{-b_j}\over (b_j)_q}
\prod_{\scriptstyle i=1\atop\scriptstyle i\neq j}^{m+1}
{(b_i-b_j+1)_q\over (b_i-b_j)_q}
\ -\
{q^{-b_{m+1}}\over (b_{m+1})_q}
\left(
\prod_{i=1}^{m}
{(b_i-b_{m+1}+1)_q\over (b_i-b_{m+1})_q}
\right)
\\[2mm]
\nonumber
&=& -
\sum_{j=1}^{m+1} {q^{-b_j}\over (b_j)_q}
\prod_{\scriptstyle i=1\atop\scriptstyle i\neq j}^{m+1}
{(b_i-b_j+1)_q\over (b_i-b_j)_q}\,
\ea
which proves formula (\ref{a1}) in the case $k=m+1$.

In order to prove eqs.(\ref{a2}), (\ref{a3})
we rewrite eq.(\ref{a1}), inverting the parameter
$q\rightarrow q^{-1}$
\be
\lb{a6}
q^{-k}\ -\ \prod_{i=1}^k {(b_i+1)_q\over (b_i)_q}\
=\ -\,
\sum_{j=1}^k {q^{b_j}\over (b_j)_q}
\prod_{\scriptstyle i=1\atop\scriptstyle i\neq j}^{k}
{(b_i-b_j+1)_q\over (b_i-b_j)_q}\, ,
\ee
and form a linear combination
$\left((\mbox{\ref{a1}})\cdot q^x -
(\mbox{\ref{a6}})\cdot q^{-x}\right)/(q-q^{-1})$, where
$x$ takes on integer values.
The resulting equality reads
\be
\lb{a7}
(k+x)_q\ -\ (x)_q\prod_{i=1}^k {(b_i+1)_q\over
(b_i)_q}\ =\
\sum_{j=1}^k {(b_j-x)_q\over (b_j)_q}
\prod_{\scriptstyle i=1\atop\scriptstyle i\neq j}^{k}
{(b_i-b_j+1)_q\over (b_i-b_j)_q}\, .
\ee
The relations (\ref{a2}) and (\ref{a3}) are particular cases of
the relation (\ref{a7}) for $x=0$ and $x=-k$, respectively.

\addtocontents{toc}{\contentsline {section}{\numberline {} References}{\pageref{refer}}}

\end{document}